\documentclass[12pt]{ucthesis}

\usepackage{latexsym}

\def \R {{I \! \! R}} 

\newtheorem{conjecture}{Conjecture}
\newtheorem{theorem}{Theorem}[chapter]
\newtheorem{lemma}[theorem]{Lemma}

\newtheorem{corollary}[theorem]{Corollary}

\input {epsf}

\def\dsp{\def\baselinestretch{2.0}\large\normalsize}
\def\bsp{\def\baselinestretch{1.5}\large\normalsize}
\dsp

\begin{document}


\title{Planar Soap Bubbles}
\author{RICHARD PAUL DEVEREAUX VAUGHN}
\degreeyear{1998}
\degree{DOCTOR OF PHILOSOPHY}
\chair{Professor Joel Hass}
\othermembers{Professor Dmitry Fuchs\\
Associate Professor Abigail Thompson}
\numberofmembers{3}
\prevdegrees{B.S. (University of San Francisco) 1990\\
M.A. (University of California, Davis) 1994}
\field{Mathematics}
\campus{Davis}

%
\begin{abstract}
\thispagestyle{empty}
The generalized soap bubble problem seeks the least perimeter way to enclose 
and separate $n$ given volumes in $\R^m$. We study the possible 
configurations for perimeter minimizing bubble
complexes enclosing more than two regions. We prove that perimeter minimizing 
planar bubble complexes with equal pressure regions and without empty chambers 
must have connected regions. As a consequence, we show that the least 
perimeter planar graph that encloses and separates three equal areas in 
$R^2$ using convex cells and without empty chambers is a 
``standard triple bubble'' with connected regions.

\vspace{.5in}

\abstractsignature
\end{abstract}

\begin{frontmatter}

\maketitle

\begin{abstract}
The generalized soap bubble problem seeks the least perimeter way to enclose 
and separate $n$ given volumes in $\R^m$. We study the possible 
configurations for perimeter minimizing bubble
complexes enclosing more than two regions. We prove that perimeter minimizing 
planar bubble complexes with equal pressure regions and without empty chambers 
must have connected regions. As a consequence, we show that the least 
perimeter planar graph that encloses and separates three equal areas in 
$R^2$ using convex cells and without empty chambers is a 
``standard triple bubble'' with connected regions.

\vspace{.5in}

\abstractsignature
\end{abstract}



\tableofcontents
\listoffigures


\begin{acknowledgements}

I want to express my deepest gratitude to my advisor and mentor, Joel Hass, for
all of his support, encouragement, and guidance. He has been extremely patient
and never complained
when I'd knock and say ``Just one quick question\ldots ''

I also want to thank the many other professors at UC Davis who have 
contributed to my education. I've been fortunate to work with a first rate 
group of topologists and geometers. I especially appreciate the contributions 
of Dmitry Fuchs for teaching me Topology and Differential Geometry and 
Abigail Thompson for many significant conversations. I also wish to thank 
them both for taking the time to read this dissertation.

I'll always be fond of the Graduate Coordinator Kathy LaGuisa. She made 
me feel welcome from the very beginning and has taken care of me ever since.

I am grateful to the many colleagues at UC Davis that have provided friendship 
and emotional support. My thanks to Maike Meyer, Ramin 
Naimi, Tanya Seph, and Michelle Stocking for all 
of their suggestions and comraderie. Curtis Feist will always be a special 
friend for his insightful comments and advice. He has contributed 
greatly to my development both personally and professionally.

Finally, I want to acknowledge my wife, Paula, who is my best friend and source
of support. She has been very understanding and patient. I am most grateful 
to Paula, her family, and my family for their 
unending love and encouragement.
\end{acknowledgements}

\end{frontmatter}

\pagestyle{myheadings}

\headheight 12pt
\headsep 25pt

\chapter{Introduction}
\thispagestyle{myheadings}

The classic \textit{isoperimetric problem} is to find the 
largest amount of area that 
can be enclosed using a simple closed curve of fixed length. The 
answer is of course 
a circle, although the proof is more difficult than some may realize.
The ancient Greeks knew the problem 
and the solution. In fact, Pappus records that 
Zenodorus found the solution first \cite{Polya}. The first 
mathematical proof, however, is
credited to Steiner \cite{Steiner} in the 19th century. He proved that if
a solution exists, then it must be a circle. Caratheodory\footnote{Blaschke 
\cite{Cara} credits Edler, Carath\'eodory and Study with existence results. 
Bandle \cite{Bandle} claims Carath\'eodory was first. Schmidt and 
Weierstrauss completed the three dimensional analogue.} 
completed the proof by showing that a solution does exist. We refer 
to the excellent prefaces in P\'{o}lya and Szeg\"{o}
\cite{Polya} and in Bandle \cite{Bandle} for more historical details.

The solution to the classic isoperimetric problem is stated by the 
isoperimetric inequality 
\[ A \leq \frac{1}{4\pi} (L)^2, \] where 
$L$ is the length of any simple closed curve in the plane and $A$ is the 
enclosed area. Equality holds if and only if the simple closed curve is a 
circle.

\begin{figure}[hbt]
\centerline {\epsfbox {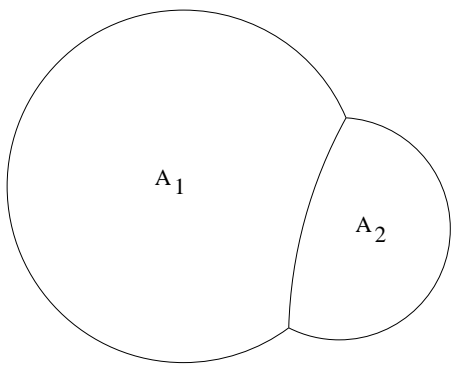}}
\caption {A standard double bubble enclosing and separating areas $A_1$ and $A_2$}
\label{double}
\end{figure}

A related problem is to find the simple closed curve with least perimeter 
that encloses a given area. 
For a single area in $\R^2$ or a single volume in $\R^m$ the problem is 
equivalent to the classic isoperimetric problem. What if we wanted, however, 
to enclose and separate two different prescribed areas in the plane? The 
solution is not two disjoint circles as we can use less perimeter by letting 
the two different areas share some of their perimeter. 
The solution \cite{Foisy} is a \textit{standard double bubble} 
consisting of two chambers enclosed by three arcs of circles all meeting at 
angles of $\frac{2\pi}{3}$. (See Figure~\ref{double}.)

The problem that seeks the least 
perimeter closed curve, surface, or hypersurface that encloses and 
separates $n$ given volumes in $\R^m$ is called 
the \textit{generalized soap bubble problem}. 
The name comes from the fact that soap 
bubbles minimize surface tension for the fixed volumes of 
air enclosed.

Another related problem is to separate or tile $\R^m$ into 
equal volume pieces 
as efficiently as possible. In $\R^2$, the \textit{honeycomb conjecture} 
states that a 
regular hexagonal tiling is the least perimeter tiling that
separates the plane into unit area pieces, although it is not immediately 
clear what is meant by least perimeter in this infinite region context. 
One interpretation is to consider this problem as the limit of the $n$ soap 
bubble problem in $\R^2$ as $n \rightarrow \infty$.

Only recently has progress been made on even the smallest cases of the 
generalized soap bubble problem. In 1976, Almgren 
\cite{Alm} proved the existence and regularity of a solution to the 
generalized soap bubble problem in dimensions bigger than 2. Taylor\cite{Tay} 
improved 
this result for dimension 3 in the same year. In 1992, Morgan \cite{Morgan} 
proved existence and regularity of a solution in dimension 2.

The double bubble problem in $\R^2$ (the least perimeter embedded planar graph 
that encloses and separates two given areas in the plane) was 
solved in 1994 by a group of 
undergraduates led by Frank Morgan \cite{Foisy}. The $n$ bubble problem in 
$\R^2$ for $n>2$, i.e. the least perimeter graph that encloses and 
separates $n$ 
given areas in $\R^2$ is the focus of this paper. 
The general $n$ bubble conjecture is that the least perimeter solution will 
always have connected regions. In 
particular, we examine three regions and the corresponding triple bubble 
conjecture.

\bigskip

\noindent \textbf{Conjecture~\ref{conj1} } \textit{The least 
perimeter planar graph that encloses 
and separates three finite areas $A_1$, $A_2$, and $A_3$ is a regular triple 
bubble complex with four vertices, six edges, and three connected regions.}

\bigskip

\begin{figure}[hbt]
\centerline {\epsfbox {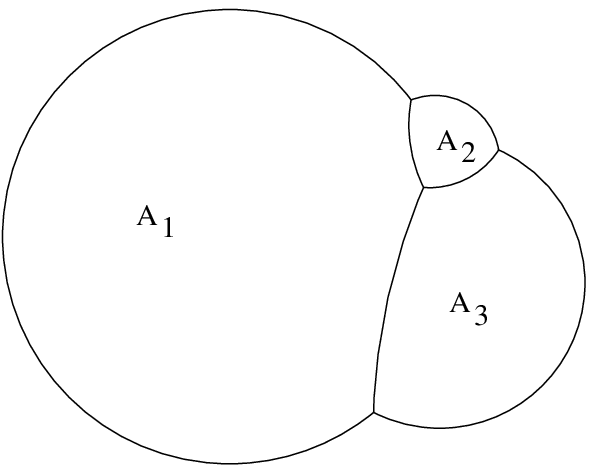}}
\caption{The conjectured minimum graph enclosing areas $A_1$, $A_2$, and $A_3$}
\end{figure}

Although much of this work was completed in the context of the triple bubble 
problem, many results are true of 
planar bubbles in general. We will specify when a result is valid only 
for three regions.

We begin Chapter 2 by looking at 
some general considerations about perimeter minimizing planar bubble 
complexes. After some definitions, we examine Morgan's existence and 
regularity theorem. We specifically look into the length $L(A_1,A_2,A_3)$ 
needed to
enclose three areas $A_1$, $A_2$, and $A_3$. We prove

\bigskip

\noindent \textbf{Lemma~\ref{continuous} } \textit{The length function 
$L(A_1,A_2,A_3)$ is continuous.} 

\bigskip

In an attempt to decrease perimeter, it is sometimes convenient to increase
the area enclosed by a region. In Conjecture~\ref{conj2}, we modify the 
triple bubble conjecture to allow these area increases.

\bigskip

\noindent \textbf{Conjecture~\ref{conj2} } \textit{Given three finite positive 
real numbers $A_1$, $A_2$, and $A_3$, the least perimeter graph that encloses 
and 
separates three finite areas $B_1$, $B_2$, and $B_3$ such that $B_i \geq A_i$ 
for all $i$, is a standard triple bubble.}

\bigskip

Using Lemma~\ref{continuous}, we prove that the two conjectures are equivalent
in

\bigskip

\noindent \textbf{Theorem~\ref{triv} } \textit{Let three positive areas 
$A_1$, $A_2$, and $A_3$ be given. There exists a least perimeter 
triple bubble 
complex $\mathcal{B}$ that encloses and separates areas 
$B_1$, $B_2$, and $B_3$ with $B_i \geq A_i$. That is, if $\mathcal{C}$ is any 
other complex enclosing areas $C_1$, $C_2$, and $C_3$ with  
$C_i \geq A_i$, then it must use at least as much perimeter $(\ell (\mathcal{B}) \leq \ell (\mathcal{C}))$. 
Furthermore, if the minimizer $\mathcal{B}$ has 
connected regions (a standard triple bubble), then it must enclose the given 
areas $A_1$, $A_2$, and $A_3$ (i.e. $B_i=A_i$ for all i). }

\bigskip

In Chapter 3, we will look at some restrictions on the shape of
 perimeter minimizing planar bubble complexes. In particular, we prove that 
connected portions of regions must have 
more than two sides, three sided pieces are 
determined by the curvature of their edges, and certain other connected pieces 
cannot touch the exterior more than once. 

In Chapter 4, we examine some additional restrictions imposed on bubble 
complexes with equal pressure regions. We show that the number of edges that
bound a portion of a region is limited to three, four, five or 
six edges. Then we show
how a few parameters control the shapes of 3-gons, 4-gons and 5-gons.
In addition, we examine the way in which these $n$-gons can meet in a perimeter
minimizing complex. For example,
when a 3-gon shares an edge with another 3-gon, they are both adjacent to 
another 3-gon, thus creating a standard triple bubble component. 

In Chapter 5, we use the results from Chapters 3 and 4 to solve 
a restricted case of the planar 
triple bubble problem. Specifically, we prove

\bigskip

\noindent \textbf{Theorem~\ref{main} } \textit{A perimeter 
minimizing triple bubble 
complex with equal pressure regions and no empty 
chambers must be a standard triple bubble. In particular, it has connected regions.}

\bigskip

Since there exists a triple bubble complex enclosing 
equal areas with connected, equal pressure regions, we have a partial solution
to the triple bubble conjecture.

\bigskip

\noindent \textbf{Corollary~\ref{cormain} } \textit{The least perimeter 
graph that 
encloses
and separates three equal areas $A_1$, $A_2$, and $A_3$ without empty chambers
using equal pressure regions is a standard triple bubble.} 

\bigskip

The restriction to equal pressure regions is analogous to the work of 
L. Fejes T\'{o}th \cite{Toth} in the 1940's.
T\'{o}th proved that the hexagonal honeycomb is the least perimeter\footnote{ Least perimeter in this tiling context means that a limiting 
perimeter to area ratio is minimized.} tiling of 
the plane with equal area polygonal cells. Since the pressure 
difference between two regions is measured by the curvature of the connecting 
edge, the restriction to polygonal cells is similar to an equal pressure 
restriction. We note this similarity in Corollary~\ref{cortriv}. 

\bigskip

\noindent \textbf{Corollary~\ref{cortriv} } \textit{The least perimeter graph 
that encloses and separates three equal areas with
convex cells and without empty chambers is a standard triple bubble.}

\bigskip

In another corollary to Theorem~\ref{main}, we 
show a bubble complex cannot be a solution to the triple bubble problem if it 
is close to a regular complex with equal pressure disconnected regions:

\noindent \textbf{Corollary~\ref{cormain2} } \textit{Suppose 
$\{\mathcal{A}_i\}$ is a sequence of regular triple bubble complexes 
that converges in length and area to a triple bubble 
complex $\mathcal{A}$. If $\mathcal{A}$ does 
not have any empty chambers, has equal pressure regions and is not a 
standard triple bubble, then there exists an $N$ such that for any $i>N$,
$\mathcal{A}_i$ is not a perimeter minimizer for the areas it encloses.}

\bigskip

In Chapter 6, we extend the arguments used in Theorem~\ref{main} to 
$n$ bubble complexes ($n>3$) with disconnected, equal pressure regions.
We point out several complexes that cannot be solutions to 
any $n$ bubble problem. In particular, we severely restrict the possible 
configurations for a perimeter minimizing bubble complex that
encloses and separates four prescribed areas with equal pressure regions.

All of our results agree with the general soap bubble conjecture. That is, we 
have not yet found a perimeter minimizing bubble complex 
with disconnected regions. 


\chapter{Definitions and Preliminaries}
\thispagestyle{myheadings}

A graph is called \textit{finite} if it has a finite number of vertices. 
We will consider only finite, planar graphs such that every vertex has 
degree at least three. 
An embedded planar graph \textit{encloses}  
areas  $A_1$, $A_2$, \ldots , $A_n$ in $\R^2$ 
if it separates the plane
into $n+1$ regions (not necessarily connected), $n$ of which contain the 
finite areas $A_1, A_2, \ldots , A_n$ respectively. Region $n+1$ is called 
the \textit{exterior} region and contains infinite area. 
Any non-exterior region is called an \textit{interior} region.
Define an \textit{n bubble complex} to be an embedded planar graph that 
encloses some $n$ positive areas.

Given $n$ positive real numbers $A_1, A_2, \ldots , A_n$, the generalized 
soap bubble problem in $\R^2$ is to find the least perimeter $n$ bubble complex
that encloses those areas.

Define a \textit{half variation} of an $n$ bubble complex $\mathcal{B}$ to 
be a continuous
family $\{ ~ \mathcal{B}_t \mid ~ t\in [0, \epsilon)\}$ of 
$n$ bubble complexes such 
that $\mathcal{B}_0=\mathcal{B}$. Let 
$\ell(\mathcal{B})$ be the function that returns the 
length of a bubble complex $\mathcal{B}$.

Suppose $\mathcal{A}$ is an $n$ bubble complex enclosing areas 
$A_1$, $A_2$, \ldots , $A_n$. If there exists a half-variation of $\mathcal{A}$
such that the areas enclosed by each $\mathcal{A}_t$ is the same as the 
areas enclosed by $\mathcal{A}$ 
and ${\frac{d\ell(\mathcal{A}_t)}{dt}\mid}_{t=0} < 0$, then $\mathcal{A}$ 
is not the least 
perimeter way to enclose $A_1$, $A_2$, \ldots , $A_n$. The half-variation
defines a deformation of $\mathcal{A}$ that preserves area and 
yet decreases perimeter.

Frank Morgan used variational arguments to prove that for any $n$ areas, 
a perimeter 
minimizer exists and must satisfy certain regularity conditions.
\begin {theorem} (Morgan\cite{Morgan}) For any positive real areas 
$A_1,A_2,\ldots,A_n$, there exists a perimeter-minimizing embedded graph 
that encloses those areas in $\R^2$. This least perimeter graph
must satisfy the following conditions:

\begin{enumerate}

\item The graph consists of a finite number of vertices, edges, and faces;

\item Edges have constant curvature (arcs of circles or line segments);

\item Vertices are trivalent;

\item Edges meet at angles of ${2\pi} \over 3$;

\item Curves separating a specific pair of regions have the same 
curvature; and

\item Any half-variation $\mathcal{A}_t$ that preserves area must not 
initially decrease length. That is, there does not exist a half-variation 
$\mathcal{A}_t$ such that 
\[ {\frac{d\ell(\mathcal{A}_t)}{dt}\mid}_{t=0} < 0.\]

\end{enumerate}
\label{regularity}
\end{theorem}  

\bigskip

We will call these six conditions the \textit{regularity conditions}.

We define a \textit{regular n bubble complex} to be an $n$ bubble 
complex that satisfies all six of the regularity conditions.
A connected portion of any interior region (i.e. a face of the embedded 
graph enclosing a piece of one of the given areas) 
will be called an \textit{n-gon}, where $n$ is the 
number of edges that enclose the connected 
piece. Edges that
separate a part of the exterior region from another region will be called 
\textit{exterior edges}. All non-exterior edges will be called 
\textit{interior edges}. 

Suppose we have a regular $n$ bubble complex $\mathcal{A}$ that encloses areas 
$A_1,A_2,\ldots ,A_n$. If we can find a
non-regular $n$ bubble complex $\mathcal{B}$ that contains 
the same areas such that 
$\ell (\mathcal{A}) \geq \ell (\mathcal{B})$, then $\mathcal{A}$ is not the 
least perimeter way to enclose $A_1,A_2,\ldots,A_n$. 
Non-regular bubble complexes yield half-variations that preserve area yet 
decrease perimeter. In other words, there exists a complex enclosing the 
same areas with length strictly less than the length of $\mathcal{B}$ and 
therefore less than the length of 
$\mathcal{A}$. This argument will be used regularly in the proof of 
Theorem~\ref{main} and throughout Chapter 6.

As stated above, Morgan established the existence of a minimum perimeter graph
that encloses
and separates any three areas. Let $L(A_1, A_2, A_3)$ be the function that
gives the minimum length needed to enclose and separate areas $A_1$, 
$A_2$, and 
$A_3$. In Lemma~\ref{continuous} we show that this length function is continuous. First, however, we prove a lemma we will need.

\begin{lemma}
\label{rescale}
For any $A_1 > 0$, $A_2 > 0$, $A_3 > 0$ and any $\epsilon > 0$, 
there exists a $\delta > 0$ such that if 
$\left| x-A_{1} \right| \leq \delta$, $\left| y-A_{2} \right| \leq \delta$, 
and $\left| z-A_{3} \right| \leq \delta$, then rescaling any 
complex enclosing $x$, $y$, and $z$ to 
get a complex enclosing $x^{\prime}$, $y^{\prime}$, and $z^{\prime}$ where 
$x^{\prime} =A_1$ or $y^{\prime} =A_2$ or $z^{\prime} =A_3$ 
will result in at most an 
$\epsilon$ change in the areas $A_1$, $A_2$ and $A_3$, i.e. 
$\left| A_{1}-x^{\prime} \right| \leq \epsilon$, 
$\left| A_{2}-y^{\prime} \right| \leq \epsilon$, and 
$\left| A_{3}-z^{\prime} \right| \leq \epsilon$.
\end{lemma}
 
\noindent \textbf{Proof:}
This is simply a consequence of the continuity of rescaling. To be precise, 
suppose $A_1$, $A_2$, $A_3$, and $\epsilon$ are given. 
By the continuity of $\frac{xA_2}{yA_1}$ for $y\neq 0$, there exists a 
$\delta_1 > 0$ such that for 
$\left| x-A_{1} \right| \leq \delta_1$ and 
$\left| y-A_{2} \right| \leq \delta_1$, we 
have $\left| 1-\frac{xA_2}{yA_1} \right| \leq \frac{\epsilon}{A_1}$.
Similarly, there exist $\delta_2, \ldots, \delta_6 > 0$ such that
\[ \left| x-A_{1} \right| \leq \delta_2 \; \mathrm{and} \;  \left| z-A_{3} \right| \leq \delta_2 \; \Longrightarrow \; \left| 1-\frac{xA_3}{zA_1} \right| \leq \frac{\epsilon}{A_1}, \]
\[ \left| x-A_{1} \right| \leq \delta_3 \; \mathrm{and} \;  \left| y-A_{2} \right| \leq \delta_3 \; \Longrightarrow \; \left| 1-\frac{yA_1}{xA_2} \right| \leq \frac{\epsilon}{A_2}, \]
\[ \left| x-A_{1} \right| \leq \delta_4 \; \mathrm{and} \;  \left| z-A_{3} \right| \leq \delta_4 \; \Longrightarrow \; \left| 1-\frac{zA_1}{xA_3} \right| \leq \frac{\epsilon}{A_3}, \]
\[ \left| y-A_{2} \right| \leq \delta_5 \; \mathrm{and} \;  \left| z-A_{3} \right| \leq \delta_5 \; \Longrightarrow \; \left| 1-\frac{yA_3}{zA_2} \right| \leq \frac{\epsilon}{A_2}, \; \mathrm{and} \]
\[ \left| y-A_{2} \right| \leq \delta_6 \; \mathrm{and} \;  \left| z-A_{3} \right| \leq \delta_6 \; \Longrightarrow \; \left| 1-\frac{zA_2}{yA_3} \right| \leq \frac{\epsilon}{A_3}. \]

Let $\delta=\mathrm{min}\left( \delta_{1},\delta_{2},\delta_{3},\delta_{4},\delta_{5},\delta_{6}\right) $.

Suppose $\mathcal{A}$ is a complex enclosing $x$, $y$, and $z$ such that 
$\left| x-A_{1} \right| \leq \delta$, $\left| y-A_{2} \right| \leq \delta$, 
and $\left| z-A_{3} \right| \leq \delta$. 

Rescale $\mathcal{A}$ to get a complex $\mathcal{B}$ enclosing $x^{\prime}$, 
$y^{\prime}$, and $z^{\prime}$.

If $x^{\prime}=A_1$, then the complex was scaled 
by a factor of $\frac{A_1}{x}$.
So, $y^{\prime} = y \left( \frac{A_1}{x} \right)$ and 
$z^{\prime} = z \left( \frac{A_1}{x} \right)$.
The first inequality is trivial $\left| A_{1}-x \right| =0<\epsilon$. But, we also
get the inequalities 
\[ \left| A_2 - y^{\prime} \right| = \left| A_{2}\left( 1-\frac{yA_1}{xA_2}\right)  \right| \leq A_2\left( \frac{\epsilon}{A_2}\right)  = \epsilon \]
and
\[ \left| A_3 - z^{\prime} \right| = \left| A_{3}\left( 1-\frac{zA_1}{xA_3}\right) \right| \leq A_3\left( \frac{\epsilon}{A_3}\right)  = \epsilon. \]

Similarly, if $y^{\prime}=A_2$, the complex was scaled by a factor of $\frac{A_2}{y}$ and we get the inequalities
\[ \left| A_1 - x^{\prime} \right| = \left| A_{1}\left( 1-\frac{xA_2}{yA_1}\right)  \right| \leq A_1\left( \frac{\epsilon}{A_1}\right)  = \epsilon, \]
\[ \left| A_{2}-y \right| =0<\epsilon \]
and
\[ \left| A_3 - z^{\prime} \right| = \left| A_{3}\left( 1-\frac{zA_2}{yA_3}\right)  \right| \leq A_3\left( \frac{\epsilon}{A_3}\right)  = \epsilon. \]

Finally, if $z^{\prime}=A_3$, then we scaled by a factor of $\frac{A_3}{z}$ 
and we get the inequalities
\[ \left| A_1 - x^{\prime} \right| = \left| A_{1}\left( 1-\frac{xA_3}{zA_1}\right)  \right| \leq A_1\left( \frac{\epsilon}{A_1}\right)  = \epsilon, \]
\[ \left| A_2 - y^{\prime} \right| = \left| A_{2}\left( 1-\frac{yA_3}{zA_2}\right)  \right| \leq A_2\left( \frac{\epsilon}{A_2}\right)  = \epsilon, \]
and
\[ \left|A_{3}-z \right| =0<\epsilon. \]
$\Box$

\bigskip

\begin{lemma}
\label{continuous}
The length function $L(A_1, A_2, A_3)$ is continuous for all $A_i > 0$.
\end{lemma}

\noindent \textbf{Proof:}
Let $A_1$, $A_2$, $A_3$, and $\epsilon >0$ be given.

Let $\delta_1 = \frac{\epsilon^2}{36\pi}$. Let $\delta_2$ be the delta 
needed in Lemma~\ref{rescale} for $\frac{\epsilon^2}{16\pi}$. That is, if
$\left| x-A_{1} \right| \leq \delta_2$, $\left| y-A_{2} \right| \leq \delta_2$,
and $\left| z-A_{3} \right| \leq \delta_2$, then 
when we rescale any complex containing  $x$, $y$, and $z$ to get one 
containing $x^{\prime}$, $y^{\prime}$, and 
$z^{\prime}$ with $x^{\prime}=A_1$ or $y^{\prime}=A_2$ or $z^{\prime}=A_3$, we
get the inequalities $\left| A_1 - x \right| \leq \frac{\epsilon^2}{16\pi}$, 
$\left| A_2 - y \right| \leq \frac{\epsilon^2}{16\pi}$, and 
$\left| A_3 - z \right| \leq \frac{\epsilon^2}{16\pi}$.

Let $\delta =\mathrm{min}\{ \delta_1, \delta_2 \}$. Suppose 
$\left| A_{1} - x \right| \leq \delta$, 
$\left| A_{2} - y \right| \leq \delta$, and
$\left| A_{3} - z \right| \leq \delta$.

Let $\mathcal{A}_1$ be a complex that uses $L(A_1,A_2,A_3)$ perimeter 
to enclose areas 
$A_1$, $A_2$, and $A_3$. Similarly, let $\mathcal{A}$ be a complex enclosing 
areas
$x$, $y$, and $z$ with length $L(x,y,z)$.

\noindent \textbf{Case 1: } $x > A_{1}$, $y > A_{2}$, $z > A_{3}$, and 
$L(x,y,z) \geq L(A_{1},A_{2},A_{3})$.

\begin{figure}[hbt]
\centerline {\epsfbox {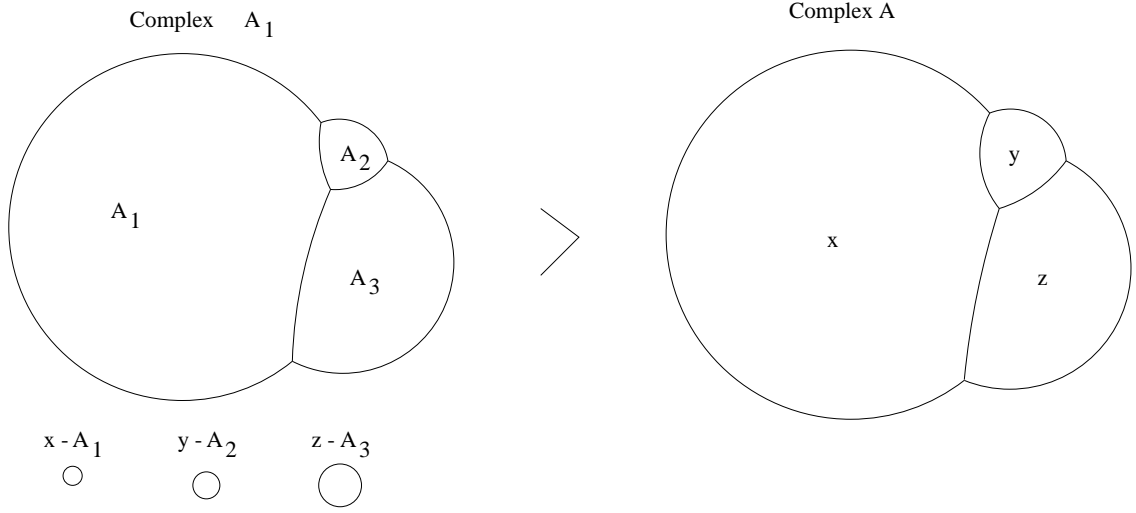}}
\caption {Adding three small circles to $\mathcal{A}_1$ gives a complex enclosing areas $x$, $y$, and $z$ with more perimeter than $\mathcal{A}$}
\label{case1}
\end{figure}

Since $x > A_{1}$, $y > A_{2}$, and $z > A_{3}$, we can 
enclose areas $x$, $y$, 
and $z$ by using $\mathcal{A}_1$ together with three disjoint 
circles containing 
areas $x-A_{1}$, $y-A_2$, and $z-A_3$ respectively. These disjoint 
circles have perimeter $2\sqrt{\pi}\sqrt{x-A_1}$, $2\sqrt{\pi}\sqrt{y-A_2}$, 
and $2\sqrt{\pi}\sqrt{z-A_3}$.
Since $L(x,y,z)$ is the minimum length, it has shorter length 
than the perimeter used by $\mathcal{A}_1$ together with the circles. (See Figure~\ref{case1}.) 

We thus get
\[ L(x,y,z) < L(A_{1},A_{2},A_{3}) + 2 \sqrt{\pi}\sqrt{x-A_{1}} + 2 \sqrt{\pi}\sqrt{y-A_{2}} + 2 \sqrt{\pi}\sqrt{z-A_{3}}. \]
So,
\[ \left| L(x,y,z) - L(A_{1},A_{2},A_{3}) \right| = L(x,y,z) - L(A_{1},A_{2},A_{3}) \]
\[ < 2\sqrt{\pi}\sqrt{x-A_{1}} + 2 \sqrt{\pi}\sqrt{y-A_{2}} + 2 \sqrt{\pi}\sqrt{z-A_{3}}\]
\[ \leq 6\sqrt{\pi}\sqrt{\delta} \]
\[ \leq 6\sqrt{\pi}\sqrt{\delta_1}\]
\[ = 6\sqrt{\pi}\sqrt{\frac{\epsilon^2}{36\pi}} = \epsilon. \]

\noindent \textbf{Case 2: } $x < A_{1}$, $y < A_{2}$, $z < A_{3}$, and 
$L(x,y,z) \leq L(A_{1},A_{2},A_{3})$.

We do exactly the same as case 1, but add little circles to $\mathcal{A}$ instead of $\mathcal{A}_1$. To be precise, we can 
enclose areas $A_1$, $A_2$, 
and $A_3$ by using $\mathcal{A}$ together with three disjoint 
circles containing 
areas $A_{1}-x$, $A_2-y$, and $A_3-z$ respectively. These disjoint 
circles have perimeter $2\sqrt{\pi}\sqrt{A_1-x}$, $2\sqrt{\pi}\sqrt{A_2-y}$, 
and $2\sqrt{\pi}\sqrt{A_3-z}$.
Since $L(A_{1},A_{2},A_{3})$ is the minimum length, it has shorter length 
than the perimeter used by $\mathcal{A}$ together with the circles. 
We thus get
\[ L(A_{1},A_{2},A_{3}) < L(x,y,z) + 2 \sqrt{\pi}\sqrt{A_{1}-x} + 2 \sqrt{\pi}\sqrt{A_{2}-y} + 2 \sqrt{\pi}\sqrt{A_{3}-z}. \]
That is,
\[ \left| L(x,y,z) - L(A_{1},A_{2},A_{3}) \right| = L(A_{1},A_{2},A_{3}) - L(x,y,z) \]
\[ < 2\sqrt{\pi}\sqrt{A_{1}-x} + 2 \sqrt{\pi}\sqrt{A_{2}-y} + 2 \sqrt{\pi}\sqrt{A_{3}-z}\]
\[ \leq 6\sqrt{\pi}\sqrt{\delta} \]
\[ \leq 6\sqrt{\pi}\sqrt{\delta_1}\]
\[ = 6\sqrt{\pi}\sqrt{\frac{\epsilon^2}{36\pi}} = \epsilon. \]

\noindent \textbf{Case 3: } $L(x,y,z)<L(A_{1},A_{2},A_{3})$ and $x>A_{1}$ or
$y>A_{2}$ or $z>A_{3}$

Scale down $\mathcal{A}$ to get a complex $\mathcal{B}$ that encloses areas
$x^{\prime}$, $y^{\prime}$, and $z^{\prime}$ such that one of areas equals 
an area from $\mathcal{A}_1$ (e.g. $x^{\prime}=A_1$) and the other two 
areas are smaller or 
equal to the remaining areas in $\mathcal{A}_1$ (e.g. $y^{\prime} \leq A_2$ and
$z^{\prime} \leq A_3$).

Let $\ell(\mathcal{B})$ be the length of $\mathcal{B}$. Since scaled minimizers are still minimizers, $\ell(\mathcal{B})=L(x^{\prime},y^{\prime},z^{\prime})$.
Also, 
since $\mathcal{B}$ was a scaled down copy of $\mathcal{A}$, we have 
$\ell(\mathcal{B}) < L(x,y,z)$.

Without loss of generality, assume that the first area is the one that is the 
same. In other words, $x^{\prime}=A_1$, $y^{\prime} \leq A_2$, and 
$z^{\prime} \leq A_3$.

\begin{figure}
\centerline {\epsfbox {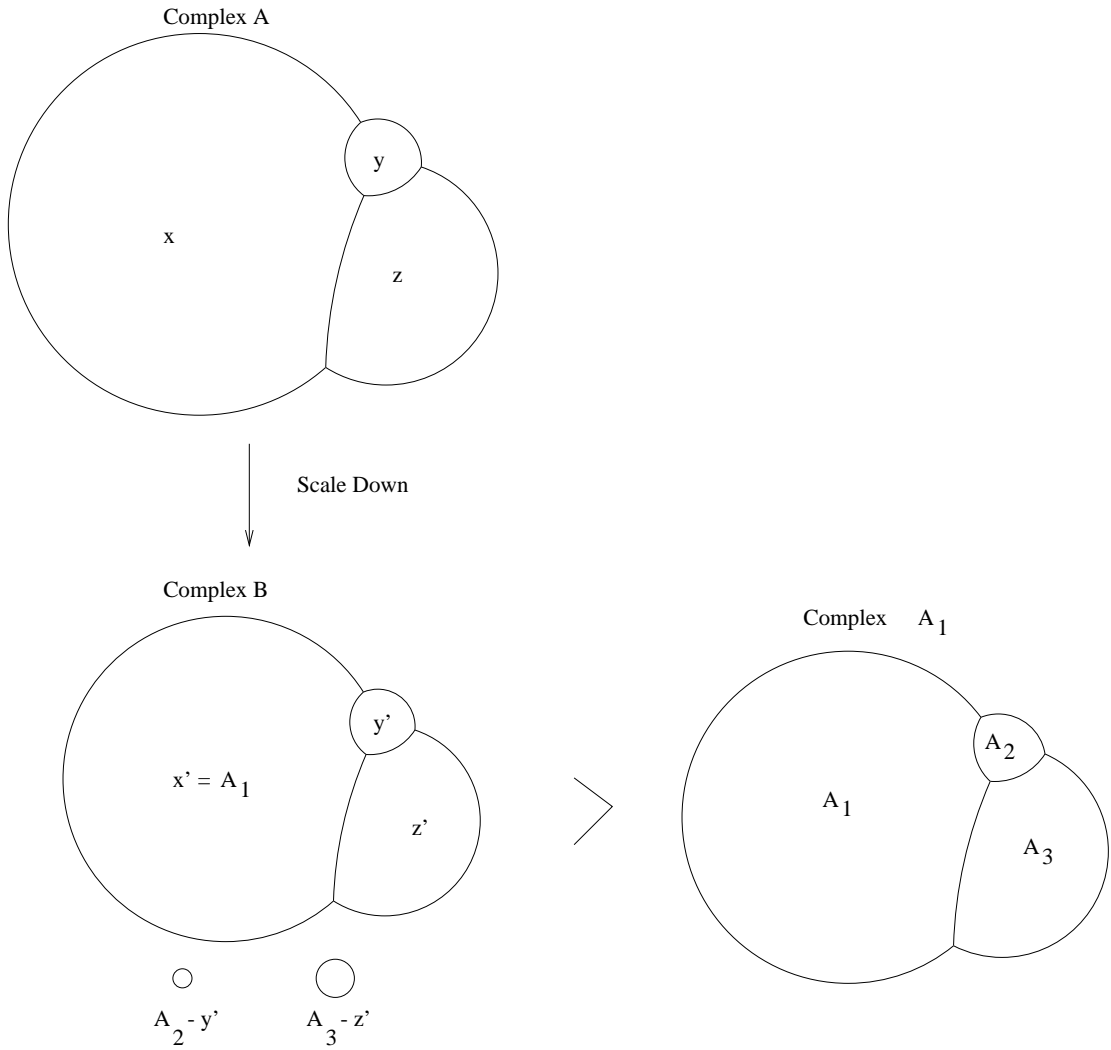}}
\caption {Scale down $\mathcal{A}$ to get $\mathcal{B}$. Then add two small circles to get a complex enclosing areas $A_1$, $A_2$, and $A_3$ with more 
perimeter than $\mathcal{A}_1$.}
\label{case3}
\end{figure}

If we add two disjoint circles of area $A_{2} - y^{\prime}$ and $A_{3} - z^{\prime}$ to the complex $\mathcal{B}$, we get a complex enclosing $A_1$, $A_2$, 
and $A_3$ again. The length of this complex is  
$\ell(\mathcal{B}) + 2\sqrt{\pi}\sqrt{A_{2} - y^{\prime}} +  2\sqrt{\pi}\sqrt{A_{3} - z^{\prime}}$ 
and must be larger than $L(A_1,A_2,A_3)$. (See Figure~\ref{case3}.) 

The 
$2\sqrt{\pi}\sqrt{A_{2} - y}$ and $2\sqrt{\pi}\sqrt{A_{3} - z}$ terms are 
the perimeter needed to 
add back in the missing area. Since we have a bound on the missing area 
(by $\delta_2$ and lemma~\ref{rescale}), we have a bound on the amount of 
perimeter needed. In fact, we get 
\[ L(A_{1},A_{2},A_{3}) < \ell(\mathcal{B}) + 2\sqrt{\pi}\sqrt{A_{2} - y^{\prime}} +  2\sqrt{\pi}\sqrt{A_{3}- z^{\prime}}\]
\[ \leq  \ell(\mathcal{B}) + 2\sqrt{\pi}\sqrt{\frac{\epsilon^2}{16\pi}} +  2\sqrt{\pi}\sqrt{\frac{\epsilon^2}{16\pi}}\]
\[ = \ell(\mathcal{B}) + \epsilon \]
In short, $L(A_{1},A_{2},A_{3}) \leq \ell(\mathcal{B}) + \epsilon$. 
Put this together with 
the previous inequalities to get
\[ \ell(\mathcal{B}) < L(x,y,z) < L(A_{1},A_{2},A_{3}) < \ell(\mathcal{B}) + \epsilon. \]
So, $\left| L(x,y,z) - L(A_{1},A_{2},A_{3}) \right| =  L(A_{1},A_{2},A_{3})-L(x,y,z) < \epsilon$.

\noindent \textbf{Case 4: } $L(x,y,z)>L(A_{1},A_{2},A_{3})$ and $x<A_{1}$ or
$y<A_{2}$ or $z<A_{3}$

This case is similar to case 3, except we scale up $\mathcal{A}$ and add little circles to 
$\mathcal{A}_1$ to get the desired inequality.

Scale up $\mathcal{A}$ to get a complex $\mathcal{B}$ that encloses areas
$x^{\prime}$, $y^{\prime}$, and $z^{\prime}$ such that one of these areas 
equals 
an area from $\mathcal{A}_1$ (e.g. $x^{\prime}=A_1$) and the other two 
areas are larger or 
equal to the remaining areas in $\mathcal{A}_1$ (e.g. $y^{\prime} \geq A_2$ and
$z^{\prime} \geq A_3$).

Let $\ell(\mathcal{B})$ be the length of $\mathcal{B}$. Since scaled minimizers are still minimizers, $\ell(\mathcal{B})=L(x^{\prime},y^{\prime},z^{\prime})$.
Also, 
since $\mathcal{B}$ was scaled up from $\mathcal{A}$, we have 
$\ell(\mathcal{B}) > L(x,y,z)$.

Without loss of generality, assume that the first area is the one that is the same. So, $x^{\prime}=A_1$, $y^{\prime} \geq A_2$, and $z^{\prime} \geq A_3$.

If we add two disjoint circles of area $y^{\prime} - A_{2}$ and 
$z^{\prime} - A_{3}$ to the complex $\mathcal{A}_1$, we get a complex 
enclosing $A_1$, $y^{\prime}$, 
and $z^{\prime}$ again. The length of this complex is  
$\ell(\mathcal{B}) + 2\sqrt{\pi}\sqrt{y^{\prime} - A_{2}} +  2\sqrt{\pi}\sqrt{z^{\prime} - A_{3}}$ 
and must be larger than $L(A_1,x^{\prime},y^{\prime})=\ell(\mathcal{B})$. 
We then get the inequality 
\[ \ell(\mathcal{B}) < L(A_{1},A_{2},A_{3}) + 2\sqrt{\pi}\sqrt{y^{\prime} - A_{2}} +  2\sqrt{\pi}\sqrt{z^{\prime} - A_{3}}\]
\[ \leq L(A_{1},A_{2},A_{3}) + 2\sqrt{\pi}\sqrt{\frac{\epsilon^2}{16\pi}} +  2\sqrt{\pi}\sqrt{\frac{\epsilon^2}{16\pi}}\]
\[ = L(A_{1},A_{2},A_{3}) + \epsilon. \]
In short, $L(A_{1},A_{2},A_{3}) > \ell(\mathcal{B}) - \epsilon$. 
Put this together with 
the previous inequalities to get
\[ \ell(\mathcal{B}) - \epsilon < L(A_{1},A_{2},A_{3}) < L(x,y,z) < \ell(\mathcal{B}) \]
So, 
\[ \left| L(x,y,z) - L(A_{1},A_{2},A_{3}) \right| =  L(x,y,z)-L(A_{1},A_{2},A_{3}) < \epsilon. \]

We have covered all the possibilities. Therefore, the length function 
$L(A_{1},A_{2},A_{3})$ is continuous. $\Box$

\bigskip

Although presented in the context of three areas, Lemma~\ref{rescale} and 
Lemma~\ref{continuous} can be easily extended to any number of regions. 
In particular, the length function $L(A_1, \ldots , A_n)$ is 
continuous for any $n$.

The triple bubble conjecture suggests a solution to the 3 bubble problem 
in $\R^2$.

\begin{conjecture}
\label{conj1}
The least 
perimeter planar graph that encloses 
and separates three finite areas $A_1$, $A_2$, and $A_3$ is a regular triple 
bubble complex with four vertices, six edges, and three connected regions.
\end{conjecture}

\begin{figure}
\centerline {\epsfbox {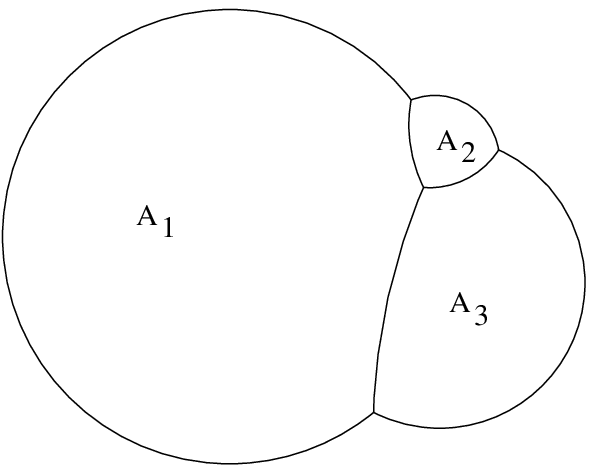}}
\caption{A standard triple bubble enclosing areas $A_1$, $A_2$, and $A_3$}
\label{tripbubb}
\end{figure}

Such a complex is 
called a \textit{standard triple bubble} and has been proven to exist 
and be unique for any three areas \cite{Montesinos}. 
(See Figure~\ref{tripbubb}.) The standard triple
bubble has also been shown to be the least perimeter way to enclose and
separate any three areas using connected regions \cite{Cox}.

An alternate version of the triple bubble conjecture allows increasing the 
areas enclosed in an attempt to minimize perimeter:

\begin{conjecture}
\label{conj2}
Given three positive 
real numbers $A_1$, $A_2$, and $A_3$, the least perimeter graph that encloses 
and 
separates three finite areas $B_1$, $B_2$, and $B_3$ such that $B_i \geq A_i$ 
for all $i$, is a standard triple bubble.
\end{conjecture}

\begin{figure}[hbt]
\centerline {\epsfbox {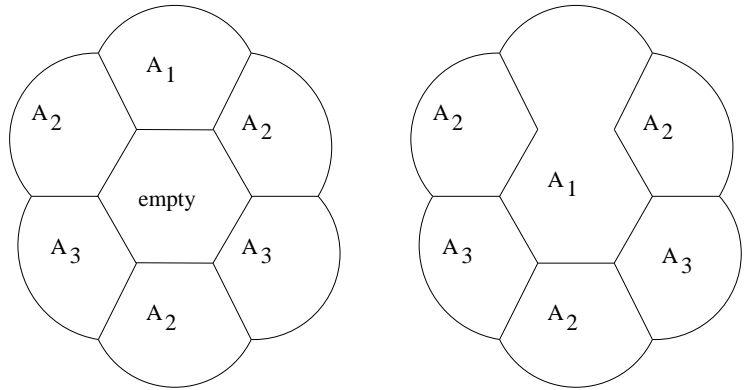}}
\caption{An empty chamber can be filled with less total length.}
\label{emcham}
\end{figure}

Conjecture~\ref{conj2} eliminates the possibility of empty 
chambers. If a complex has an empty chamber, it could be filled in with any 
one of the adjacent areas and at least one edge could be eliminated. The 
resulting complex encloses more area, but uses less perimeter. (See 
Figure~\ref{emcham}.)

Theorem~\ref{triv} establishes that the two versions of the triple bubble 
conjecture are equivalent.

\begin{theorem}
\label{triv}Let three positive areas 
$A_1$, $A_2$, and $A_3$ be given. There exists a least perimeter 
triple bubble 
complex $\mathcal{B}$ that encloses and separates areas 
$B_1$, $B_2$, and $B_3$ with $B_i \geq A_i$. That is, if $\mathcal{C}$ is any 
other complex enclosing areas $C_1$, $C_2$, and $C_3$ with  
$C_i \geq A_i$, then it must use at least as much perimeter $(\ell (\mathcal{B}) \leq \ell (\mathcal{C}))$. 
Furthermore, if the minimizer $\mathcal{B}$ has 
connected regions (a standard triple bubble), then it must enclose the given 
areas $A_1$, $A_2$, and $A_3$ (i.e. $B_i=A_i \; \; \forall i$).

\end{theorem}

\noindent \textbf{Proof:}

To optimize the length function $L$ for areas greater than or equal to 
$A_1$, $A_2$, and $A_3$, the domain we need to consider is bounded below (by 
$A_1$, $A_2$, and $A_3$) and bounded above as well. The upper bound can be 
chosen to be $B_1$, $B_2$, $B_3$ where 
$B_i=(2\pi+3)\frac{\sqrt{A_1+A_2+A_3}}{\sqrt{\pi}}$. This is the area needed 
to enclose and separate the three areas with a circle and three radii. Since
the total perimeter needed to enclose even one $B_i$ is larger than a known 
way to enclose $A_1$, $A_2$, and $A_3$, the total 
perimeter used in any attempt to enclose areas bigger than $B_1$, $B_2$, 
and $B_3$ must be larger than the minimum way to enclose $A_1$, $A_2$, and $A_3$. Since a continuous function on a 
compact set achieves its maximum and 
minimum value, there is a minimum value for the length function.

\begin{figure}[hbt]
\centerline {\epsfbox {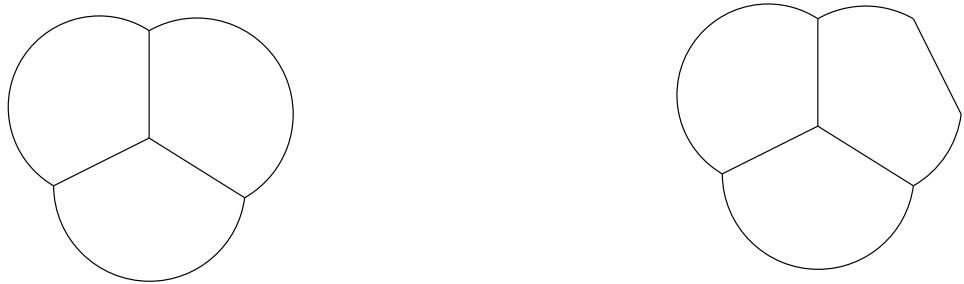}}
\caption{Cut off a little bit to save perimeter.}
\label{shrink}
\end{figure}

If the minimum is always a standard triple bubble (i.e. Conjecture 2 is 
correct), then the minimum must enclose exactly 
$A_1$, $A_2$, and $A_3$. If it encloses some $B_1$, $B_2$, and 
$B_3$ with $B_i>A_i$ for some $i$, we could reduce perimeter by replacing
a small portion of the exterior arc of $B_i$ by a straight line. (See 
Figure~\ref{shrink}.) The line can 
be chosen small enough so that the area enclosed by the region is still larger 
than $A_i$, and yet we've used less perimeter. This contradicts the assumption 
that the complex was the minimum. $\Box$

\bigskip

\begin{figure}[hbt]
\centerline {\epsfbox {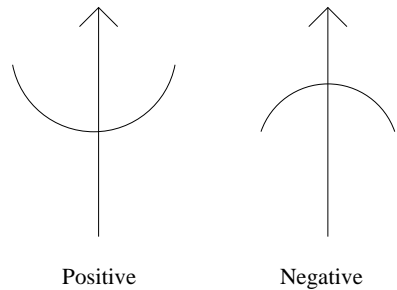}}
\caption{Sign convention for curvature}
\label{orient}
\end{figure}

We define the \textit{pressure} of a region in a regular $n$ bubble complex to 
be 0 for the exterior region. For any other region, we pick a path from the 
exterior
to that region such that the path intersects the edges of the complex
transversely in a finite number of points. The pressure is then the sum of the 
signed curvatures of the edges at these finite number of intersection points. 
We use the sign convention as shown in Figure~\ref{orient}. When exterior 
edges bulge outward (as in soap bubbles), the choice of sign
guarantees that regions adjacent to the exterior have positive pressure. It
also makes the sign of the curvature agree with the standard definition of 
curvature when the edges are given a counter-clockwise orientation.

Cox, Harrison, Hutchings, et. al. \cite{Cox} proved that for any 
closed path intersecting a regular bubble transversely, 
the sum of the signed 
curvatures along that path must be zero. In Lemma~\ref{reg}, we generalize this
result to any path that starts and ends in the same (possibly disconnected) 
region. It also guarantees that pressure is well defined.

\begin{lemma}
\label{reg}
Let $\mathcal{A}$ be a regular $n$ bubble complex. Let $\gamma$ be any path 
that intersects the edges of the complex transversely such that $\gamma$ 
starts and stops in portions of the same region (not necessarily connected).
Then, the sum of the signed curvatures of the edges crossed is zero.
\end{lemma}

\noindent \textbf{Proof:} 
Suppose that $\gamma$ goes through regions $R_1, R_2, \ldots, R_n, R_1$ and 
crosses edges with curvatures $\kappa_1, \kappa_2, \ldots, \kappa_n$ at 
points $p_1, p_2, \ldots, p_n$. Define a half-variation $\mathcal{A}_t$ 
that transfers
$t$ area from each $R_i$ to $R_{i+1}$ by adjusting each edge in a neighborhood 
about $p_i$. The initial change in length by this half-variation is just the 
sum of the signed curvatures (see e.g. Morgan\cite{M2}). That is, 
\[ \frac{d\ell (\mathcal{A}_t)}{dt}\mid_{t=0} = \sum_{i=1}^{n} \kappa_i. \]
By regularity (condition six), this sum must be greater than or equal to 0. 
If, however, the sum is greater than zero, we can traverse $\gamma$ in the 
opposite direction to get the same curvatures with 
opposite orientation. Therefore, the 
half-variation defined by adjusting area along $-\gamma$ has negative initial 
change in length
which violates regularity. Thus, the sum of the curvatures must be zero.
$\Box$

\bigskip

In Chapter~\ref{result}, we show that Conjecture 2 is true in 
the case of equal pressure regions, or that Conjecture 1 is true in 
the case of equal pressure regions with no empty chambers. 
The restriction that 
the regions have equal pressures guarantees that the inner edges (edges that 
don't touch the connected exterior region) are all line segments (0 curvature)
and the outer edges all have the same curvature. In particular, every 
$n$-gon is convex.

\chapter{Structure of Perimeter Minimizing Bubbles}
\thispagestyle{myheadings}

We begin with some observations about the possible configurations for 
perimeter minimizing bubble complexes. Theorem~\ref{regularity} guarantees
that they must be regular bubble complexes. The restrictions we discuss in this chapter are applicable to arbitrary regular bubble complexes 
enclosing any number of regions.

We first note that perimeter minimizing complexes must be connected. If a 
complex has two disconnected components, they can be pushed together until a 
vertex of degree at least four is created. This new complex violates 
regularity and therefore there exists a complex enclosing the same areas
with less perimeter.

\begin {lemma} Perimeter minimizing regular $n$ bubble complexes ($n>2$) have 
no 2-gons.
\label{lem2gon}
\end{lemma}

\begin{figure}[hbt]
\centerline {\epsfbox {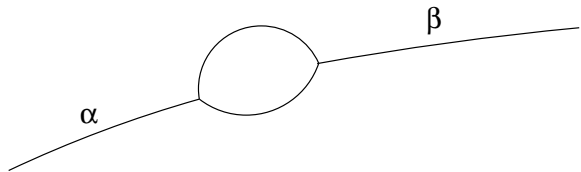}}
\caption{A 2-gon in a bubble complex}
\label{bubb1}
\end{figure}

\noindent \textbf{Proof:} Suppose there is a 2-gon. 
By regularity, every vertex must be trivalent. In 
particular, a 2-gon will have two vertices and two edges with an additional 
edge leading away from each vertex $\alpha$ and $\beta$. 
(See Figure~\ref{bubb1}.)


\begin{figure}[hbt]
\centerline{\epsfbox {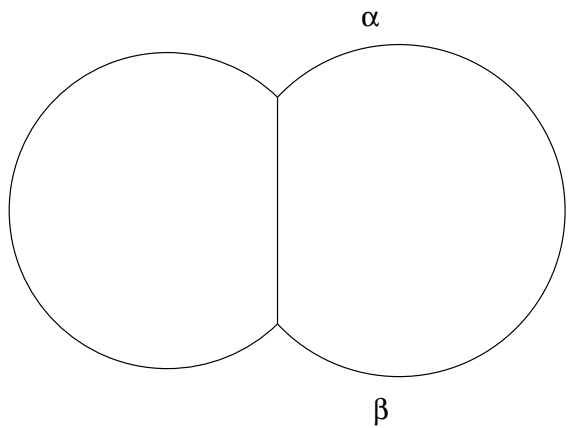}}
\caption{A disconnected double bubble region}
\label{bubb2}
\end{figure}

\noindent Case I:  edge $\alpha =$ edge $\beta$ (See Figure~\ref{bubb2}.)


The 2-gon and the adjacent 2-gon form a double bubble disconnected 
from the rest of the complex. Move this disconnected
piece until it touches another component of the bubble complex. A four valent 
vertex would be created at the point of intersection thus violating regularity.

\noindent Case II: edge $\alpha \neq$ edge $\beta$ (See Figure~\ref{bubb1}.)

The data from a single vertex is enough to completely determine a connected 
double bubble complex. In other words, if three arcs of circles meet at a 
vertex at equal angles $(\frac{2\pi}{3})$ and the sum of the signed curvatures 
of the arcs around the vertex is zero, then the arcs will extend to a 
standard double bubble complex.
All three arcs meet again at some other point and with the same angles as 
the angles at which they leave $(\frac{2\pi}{3})$.

\begin{figure}[hbt]
\centerline{\epsfbox{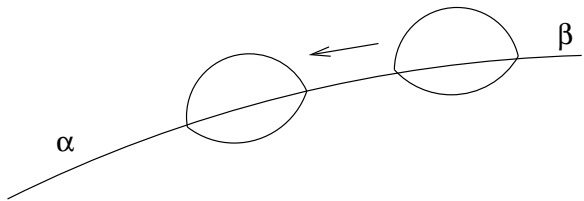}}
\caption{Slide a 2-gon along a circle}
\label{slide}
\end{figure}

By regularity, the curvatures of $\alpha$ and $\beta$ are determined by the 
curvatures of the 2-gon and therefore must be the same. Furthermore, edge 
$\alpha$ and edge $\beta$ must be arcs of the same circle since the data from 
one vertex is enough to determine the other. The 2-gon can be \textit{slid} 
along this circle without changing perimeter or area. That is, we can remove 
the 2-gon and extend edge $\alpha$ and $\beta$ to get a continuous arc of a 
circle. The 2-gon can then be reinserted anywhere along this arc. (See 
Figure~\ref{slide}.) To finish the slide move, we erase the portion of the circle inside of the new 2-gon. We then have a bubble complex enclosing equivalent
areas with exactly the same amount of perimeter. The slide move was 
introduced by the SMALL Geometry Group \cite{Foisy} when they proved the 
planar double bubble conjecture.

We continue sliding this 
2-gon until it either touches another edge or the edge $\alpha$ disappears. 
In either case, a 4-valent vertex is created and regularity is violated.
$\Box$

\bigskip
\begin{lemma}
If there exists a 3-gon with edges of curvature $\kappa_1$, $\kappa_2$, and 
$\kappa_3$, then its shape is unique (up to orientation and isometry) 
and is determined by the triangle of its vertices.
\end{lemma}

\noindent \textbf{Proof:}
When two arcs of circles of radius $r_1$ and $r_2$ meet at an angle 
of $\frac{2\pi}{3}$, the centers of the two circles are at distance
$d=\sqrt{r_1^2 + r_2^2 -r_1 r_2}$. Suppose that there exists a 3-gon
with curvatures $\kappa_1$, $\kappa_2$, and $\kappa_3$.

If two curvatures are zero, the 3-gon is determined by the curvature of
the third arc. Uniqueness is guaranteed by the Gauss-Bonnet Theorem (See Lemma 4.2).

\begin{figure}
\centerline{\epsfbox{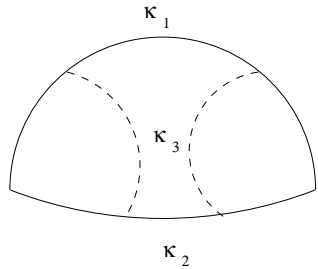}}
\caption{Two oppositely oriented 3-gons can be created from a 2-gon by adding an arc of given curvature $\kappa_3.$}
\label{uniq3gon}
\end{figure}

If none of the curvatures are zero, the radii of the respective arcs are $r_i=\frac{1}{\kappa_i}$. Consider the 2-gon 
formed by intersecting a circle of radius $r_1$ with a circle of radius $r_2$
at an angle of $\frac{2\pi}{3}$. Let $C_1$ be the circle obtained by 
extending the arc 
with curvature $\kappa_1$. Similarly, let $C_2$ denote the circle obtained
by extending the arc with curvature $\kappa_2$. To get any 3-gon with the 
same curvature edges,
we need to add a third circle of radius $r_3$ so that the angles made with 
both circles $C_1$ and $C_2$ is again $\frac{2\pi}{3}$. 
Since there is a 3-gon with these curvatures, we know it is possible.
The center of this circle must be at
distance $d_1=\sqrt{r_1^2 + r_3^2 - r_1 r_3}$ from the center of $C_1$ and
distance $d_2=\sqrt{r_2^2 + r_3^2 - r_2 r_3}$ from the center of $C_2$. 
Construct a circle of radius $d_1$ around the center of $C_1$ and a circle of 
radius $d_2$ around the center of $C_2$. The center of the third circle must lie on the intersection of these circles. The circles are not equivalent 
(since they have different centers) and intersect at least once (since there is a solution). The only other possibility is that the circles intersect twice.
If so, the two choices for the center of the third circle give the same 
intersection pattern with the 2-gon, but on 
opposite sides. (See Figure~\ref{uniq3gon}.) Generically, when an arc 
can be added to a 2-gon 
to form a 3-gon using one of the original vertices of the 2-gon, an arc of 
the same curvature can be added as well with 
opposite orientation. This creates a 3-gon with the opposite vertex 
of the 2-gon. The two different 3-gons created are equivalent 
but have opposite orientation.

If only one curvature is zero, we can build a 2-gon with a straight line and a circle of curvature $\kappa_2 \neq 0$. Consider arcs of curvature $\kappa_3$ 
leaving a point on the line segment between the vertices of our 
2-gon at an angle of $\frac{2\pi}{3}$. When such an arc meets the given arc 
of curvature $\kappa_2$, the angle made is strictly increasing 
between $0$ and $\pi$ as the 
point of departure varies from one vertex to the other. At only one point is 
the angle exactly $\frac{2\pi}{3}$.

\begin{figure}[hbt]
\centerline{\epsfbox {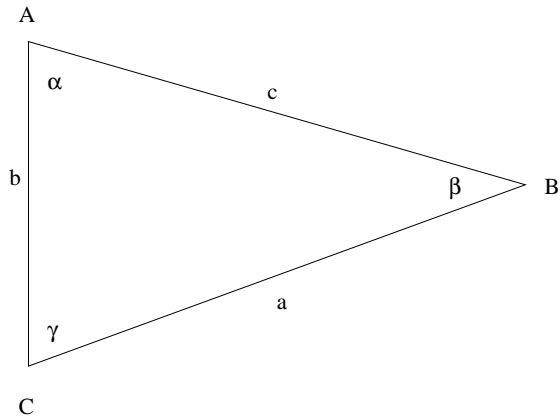}}
\caption{A generic triangle}
\label{triang}
\end{figure}

To prove that 3-gons are determined by triangles, 
we will establish a map from triangles to 3-gons and 
show that it is bijective. Suppose we have a triangle $\triangle ABC$ with 
side lengths $a$, $b$, and $c$ (opposite side from the appropriately labeled 
vertex) and angles $\alpha =\angle CAB$, $\beta=\angle ABC$, and 
$\gamma=\angle ACB$. (See Figure~\ref{triang}.)

\begin{figure}[hbt]
\centerline{\epsfbox {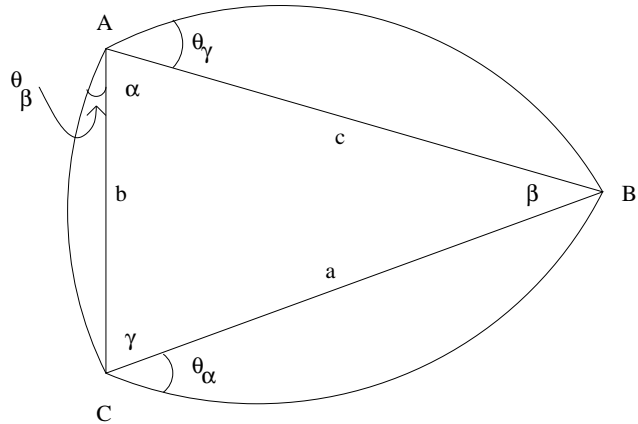}}
\caption{A triangle with arcs of circles attached}
\label{triang1}
\end{figure}

Given any angle $\theta_{\alpha}$, there is a unique arc of a circle that 
passes
through $B$ and $C$ that makes angle $\theta_{\alpha}$ with the line segment 
$BC$. In
fact, since $a$ is the length of $BC$, the curvature $\kappa_{\alpha}$ of 
the arc 
through $BC$ 
with angle $\theta_{\alpha}$ is given by the formula 
$\kappa_{\alpha} = \frac{2sin(\theta_{\alpha})}{a}$. We
consider angles exterior to the triangle to be positive.
(See Figure~\ref{triang1}.)
Similarly, angles $\theta_{\beta}$ and $\theta_{\gamma}$ uniquely determine 
arcs of circles through $AC$ and $AB$ respectively.

To get a valid 3-gon, the internal angles should all be $\frac{2\pi}{3}$. A 
3-gon must then satisfy the linear equations
\[ \theta_{\beta} + \alpha + \theta_{\gamma} = \frac{2\pi}{3}, \] 
\[ \theta_{\alpha} + \beta + \theta_{\gamma} = \frac{2\pi}{3}, \]
and
\[ \theta_{\alpha} + \gamma + \theta_{\beta} = \frac{2\pi}{3}. \]

From the triangle we also get the equation
\[ \alpha +\beta + \gamma = \pi. \]
The unique solution to these equations is 
\[ \theta_{\alpha}=\alpha - \frac{\pi}{6}\]

\[ \theta_{\beta}=\beta - \frac{\pi}{6}\]

\[ \theta_{\gamma}=\gamma - \frac{\pi}{6}.\]

Given any 3-gon, we can get a triangle by connecting the vertices of the 3-gon.
 So, the map from triangles to 3-gons is surjective. Suppose now that two 
different triangles produce the same 3-gon. Since the vertices of the 3-gon 
coincide with the vertices of the triangle that produced it, the two 
triangles must be identical. Therefore the map is also injective.
$\Box$

\bigskip

\begin{lemma} Any 4-gon or 5-gon in a perimeter minimizing regular bubble 
complex without empty chambers shares at most one edge with the exterior 
region.
\label{lemext}
\end{lemma}

\begin{figure}[hbt]
\centerline{\epsfbox {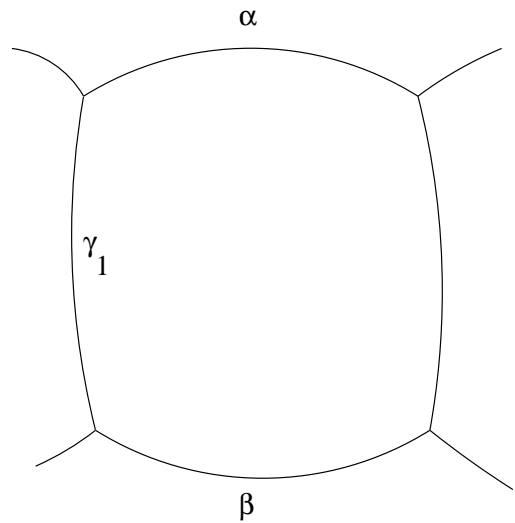}}
\caption{A region with two exterior edges}
\label{extedg}
\end{figure}

\begin{figure}[hbt]
\centerline{\epsfbox {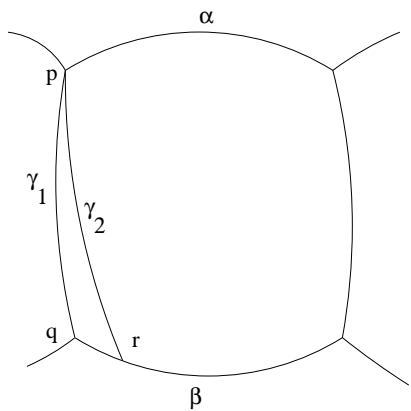}}
\caption{Cut out this wedge shaped region}
\label{wedge}
\end{figure}

\noindent \textbf{Proof:} A 4-gon or a 5-gon is distinguished by the fact that 
any pair
of edges are separated by at most one edge. Suppose a 4-gon or 5-gon shares 
two edges $\alpha$ and $\beta$ with the exterior. 
Then there is a single edge $\gamma_1$ that connects the two exterior edges. 
(See Figure~\ref{extedg}.)
Let $p$ be the vertex shared by $\gamma_1$ and $\alpha$ and 
$q$ be the vertex shared by $\gamma_1$ and $\beta$.
Pick a point $r$ at distance $\epsilon$ from $q$ on $\beta$. 
Let $\gamma_2$ be an arc of a circle or line segment from $p$ to $r$ that does 
not intersect any edges of the $n$-gon. (See Figure~\ref{wedge}.) 
If the $n$-gon is convex, the arc 
$\gamma_2$ can always be chosen to be a straight line segment.

\begin{figure}[hbt]
\centerline{\epsfbox {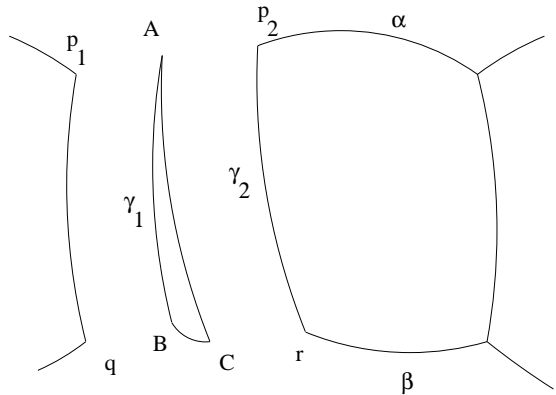}}
\caption{The cut apart complex}
\label{cutwedge}
\end{figure}

Cut out the triangular wedge formed by $\gamma_1$, $\gamma_2$, and the arc 
from $q$ to $r$. Label the corners of this triangular wedge as follows: 
$A$ for the corner that came from the point $p$, $B$ for the corner that came 
from the point $q$, and $C$ for the corner that came from the point $r$.
Then separate the remaining complex into three disjoint pieces 
by splitting $p$ into two points $p_1$ and $p_2$ such that $p_1$ is connected 
to the point $q$ and $p_2$ is connected to the point $r$. 
(See Figure~\ref{cutwedge}.) 

\begin{figure}[hbt]
\centerline{\epsfbox {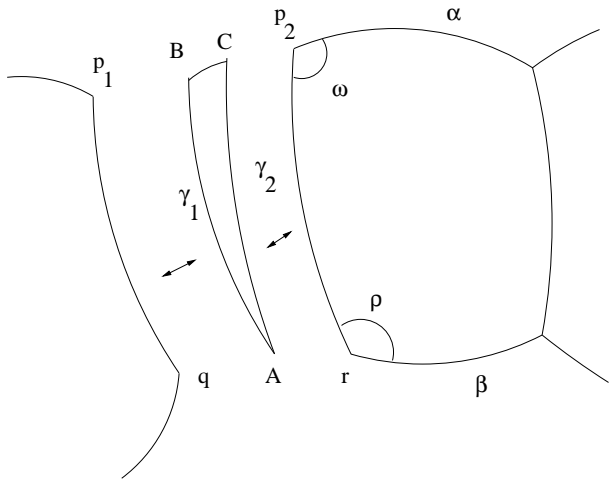}}
\caption{Flip over the wedge and glue it back in.}
\label{wedgie}
\end{figure}

Now re-attach the triangular wedge with an opposite orientation by identifying
the point $p_1$ with $B$, $p_2$ with $C$, and the points $q$, $r$, and $A$ 
with each other. (See Figure~\ref{wedgie}.)

The resulting complex has identical perimeter and encloses the same areas. The 
edge $\alpha$, however, either has a corner at $p_2$ (if angle $\omega$
differs from angle $\rho$ ) or is still smooth but longer. If the edge has a 
corner, we have a 
complex that encloses the same areas with the same perimeter that violates 
regularity. If the angles agree and no corner is created, the effect is that 
the complex has \textit{slid} distance $\epsilon$ along the edge $\alpha$. 
Continue sliding 
(i.e. repeat this procedure) until either the complex bumps into itself 
somewhere or the edge $\beta$ disappears. In either case, a complex 
that encloses the same areas with identical perimeter is created. But, this 
new complex contains a four-valent vertex and therefore violates regularity.
$\Box$

\bigskip

The same argument can be used to show that many other kinds of $n$-gon's in a 
perimeter minimizing regular bubble complex cannot touch the exterior region
more than once. All that is needed to extend the argument is the existence of
symmetric arcs $\gamma_1$ and $\gamma_2$ that intersect one exterior edge at a 
point and the other exterior edge at two points distance $\epsilon$ apart such
that the symmetric arcs do not touch any other portion of the boundary of the
$n$-gon.

In a paper regarding triple bubbles with connected regions, Cox, Harrison, 
Hutchings et. al. \cite{Cox} made an interesting remark. They found a nice 
relationship between the perimeter $L$ of a regular $n$ bubble complex, 
the areas enclosed $A_1, A_2, \ldots, A_n$, and the pressure of each region 
$p_1, p_2, \ldots, p_n$:
\begin{equation}
L=2 \sum_{i=1}^{n} p_i A_i
\label{eqhut}
\end{equation}
Using this relationship, we easily prove the following lemma:

\begin{lemma}
\label{lemhutch}
If a regular $n$ bubble complex $\mathcal{C}$ encloses areas 
$A_1, A_2, \ldots, A_n$ and 
has pressure $p_1, p_2, \ldots, p_n$ respectively, then any regular $n$ 
bubble complex $\mathcal{B}$ enclosing the same 
areas with pressures $q_i < p_i$ for all 
$i$ is not a minimizer.
\end{lemma}

\noindent \textbf{Proof:}
Let $\ell(\mathcal{B})$ be the length of $\mathcal{B}$ and $\ell(\mathcal{C})$ 
be the length of $\mathcal{C}$. Then, equation~\ref{eqhut} gives us
\[ \ell(\mathcal{B})=2\sum_{i=1}^{n} q_i A_i > 2\sum_{i=1}^{n} p_i A_i = \ell(\mathcal{C}) \]
$\mathcal{C}$ uses less perimeter and 
therefore $\mathcal{B}$ is not a perimeter minimizer. 
$\Box$

\chapter{Restrictions Imposed By Equal Pressure Regions}
\thispagestyle{myheadings}

In this chapter we consider regular $n$ bubble complexes that have 
equal pressure
interior regions. Since there is no pressure change from one interior region
to another, the curvature of the interior edges must be zero. In addition, by 
following a closed path that touches only two 
distinct exterior edges (any number
of interior edges), we get that the
curvature of the exterior edges must all be the same.

\begin {lemma} If a regular $n$ bubble complex has no empty chambers and has 
positive equal pressure regions, then $n$-gons have
at most 6 sides. In addition, $n$-gons that share an edge with the exterior 
region have at most 5 sides.
\end{lemma}

\noindent \textbf{Proof:} By regularity conditions, the internal angles of 
each $n$-gon must be 
$2\pi \over 3$. In addition, their edges have either 0 curvature 
(edges that separate two $n$-gons) or one
fixed curvature (edges that separate an $n$-gon from the exterior). 
Furthermore, since regions have positive pressure, edges separating an $n$-gon from the exterior must bulge outward (i.e. have positive curvature).

Using arcs of constant positive curvature in the Gauss-Bonnet Theorem, we get 
that
\[\sum_{i=1}^{n} \kappa_{i} l_{i} + \sum_{i=1}^{n} {\pi \over 3} = 2\pi \] 
where $n$ is the number of edges. 
Note that since each interior angle of an $n$-gon is $2\pi \over 3$, 
each exterior angle is $\pi \over 3$.
We solve for the sum of the exterior angles to get
\[ n({\pi \over 3})=2\pi-\sum_{i=1}^{n} \kappa_{i} l_{i} \leq 2\pi. \]
So, $n \leq 6$ with equality only when the edges all have 0 curvature. That is,
6-gon's are \textit{internal} since they do not share an edge with the 
exterior. 

Recall that there also can't be any 2-gons in a perimeter minimizing complex (by Lemma~\ref{lem2gon}. 
Therefore, $n$-gons that share an edge with the exterior must be 3-gons, 
4-gons, or 5-gons.
$\Box$

\bigskip

\begin{lemma}
In an $n$ bubble complex with equal pressure regions, 
there is a unique shape for a 3-gon region, a one
parameter family of possible 4-gons, and a two parameter family of possible 
5-gons (up to orientation preserving isometry). 
4-gons are determined by the length of a side adjacent to an exterior
edge. 5-gons are determined by the lengths of any two of the non-curved edges.
\label{lemstruc}
\end{lemma}

\noindent \textbf{Proof:} For a 3-gon, the Gauss-Bonnet theorem says that 
\[ \sum_{i=1}^{3} (\kappa_{i} l_i + {\pi \over 3}) = 2\pi. \]
Since two of the sides have no curvature, we get
$\kappa l_1 = \pi$ or $l_1 = {\pi \over \kappa}$ where $\kappa$ is the 
curvature
of the outside arcs and $l_1$ is the length of that curved arc in a 3-gon.
In particular, the length of the curved arc in a 3-gon with fixed 
curvature $\kappa$
is constant.

\begin{figure}[hbt]
\centerline {\epsfbox {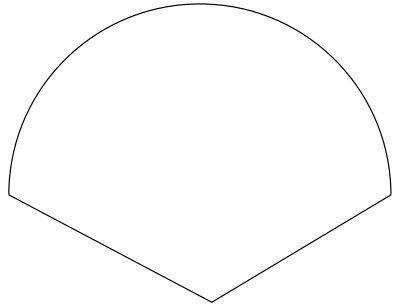}}
\caption {The unique shape for a 3-gon}
\label{3gon}
\end{figure}

 Now, two straight lines that leave the ends of such an arc at an 
angle of ${2 \pi} \over 3$ will then meet in only one point, also 
at an angle of ${2 \pi} \over 3$. This is the unique configuration for a 3-gon.
(See Figure~\ref{3gon}.)

\begin{figure}[hbt]
\centerline {\epsfbox {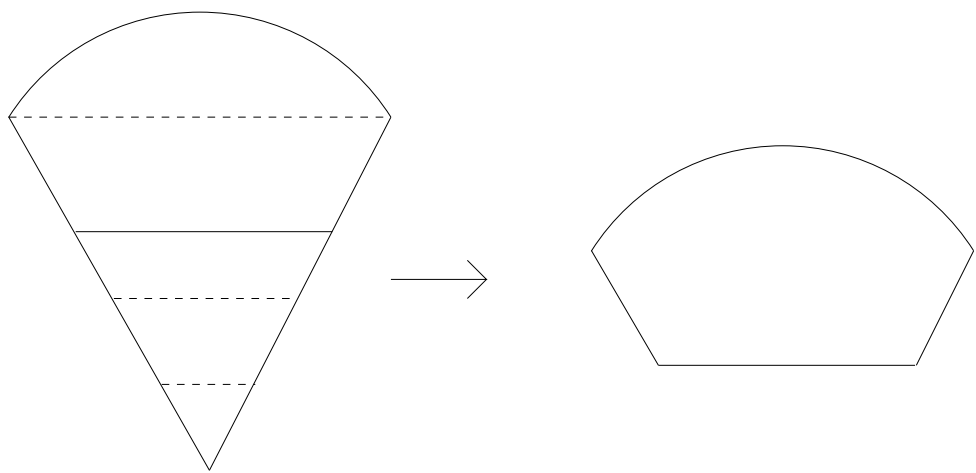}}
\caption {Cut out an equilateral triangle to determine a 4-gon}
\label{4gon}
\end{figure}

For a 4-gon, the length of the curved arc is again completely determined by 
the Gauss-Bonnet theorem. Two lines that leave the end points of this arc at 
${2 \pi} \over 3$ will meet at some point at an angle of ${\pi} \over 3$. 
Draw the line 
segment connecting the ends of the arc. The fourth side of the 4-gon (opposite 
the curved arc) must be parallel to this line segment. A choice of how far 
along a side edge to place this opposite edge completely determines the 4-gon.
This choice also corresponds to the size of equilateral triangle that is cut off of the bottom. (See Figure~\ref{4gon}.)

\begin{figure}[hbt]
\centerline {\epsfbox {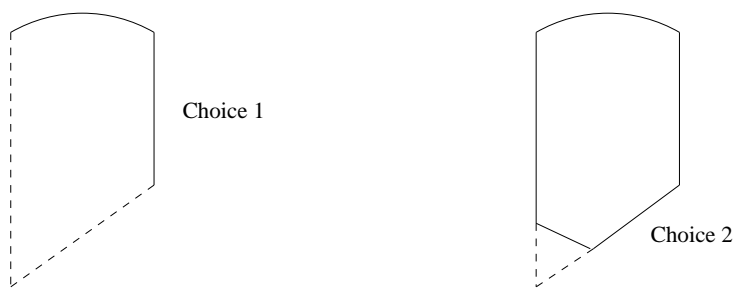}}
\caption {Two choices determine a 5 gon}
\label{5gon}
\end{figure}

For a 5-gon, the length of the arc is again fixed. This time the adjacent
edges will be parallel to each other. Choose any length for one edge, and draw
a line segment from the end of this edge such that the internal angle is 
${2 \pi} \over 3$. This edge will meet the other adjacent edge at an angle of 
${\pi} \over 3$. (See Figure~\ref{5gon}.)

Now, the vertex opposite the curved edges can be chosen to be any point on 
this line segment. Once this choice is made, the 5-gon is completely 
determined. $\Box$

\bigskip

Regular bubble complexes with equal pressure, disconnected 
regions have a lot of 
symmetry that can be used
to find non-regular complexes that enclose the same areas using equal or less
total perimeter. Since 3-gons are unique, it is convenient to consider what can
possibly be next to a 3-gon.

If a 3-gon shares an edge with another 3-gon in a minimizing complex, then 
they both share an edge with a third 3-gon. Otherwise, the region adjacent to 
both 3-gons would have two exterior edges which violates Lemma~\ref{lemext}. 
The complex is either disconnected and not a minimizer or is a standard 
triple bubble.

\begin{figure}[hbt]
\centerline {\epsfbox {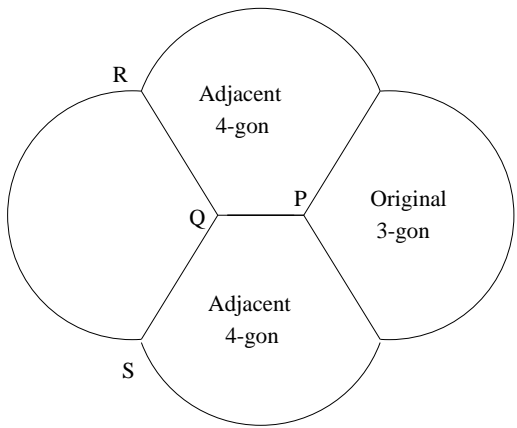}}
\caption {Two 4-gons adjacent to a 3-gon}
\label{4bubble}
\end{figure}

Suppose a 3-gon is adjacent to two 4-gons (one on each side). Let 
$P$ be the vertex shared by the 3-gon and both 4-gons. Let $Q$ be the other
interior vertex shared by the 4-gons. Let $R$ and $S$ be the vertices of the 
4-gons not shared with the original 3-gon. (See Figure~\ref{4bubble}.) 
Edges $\overline{RQ}$ and $\overline{SQ}$ separate the 4-gons from 
another $n$-gon.

This $n$-gon (\textit{opposite} the original 3-gon) must also be a 
3-gon. If not,
then the $n$-gon has two different exterior edges (One coming from vertex 
$R$ and one from vertex $S$) in violation of 
Lemma~\ref{lemext}. Since the complex is connected, there aren't any more 
regions. The complex must have just these four chambers.

Lemma~\ref{lem34} below, together with the fact that 3-gons are unique, 
proves that the complex is actually completely determined, 
i.e. it is just a scaled copy of Figure~\ref{4bubble}.

\begin{figure}[hbt]
\centerline {\epsfbox {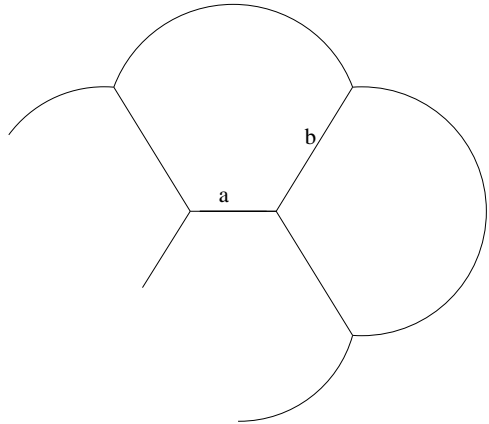}}
\caption {A 3-gon adjacent to a 4-gon}
\label{34gon}
\end{figure}

If a 3-gon shares an edge with a 4-gon, the 4-gon is unique. That is, there 
is only one possible shape for a 4-gon adjacent to a 3-gon since the one parameter has been determined for the 4-gon. (See Figure~\ref{34gon}.) Lemma~\ref{lem34} gives a nice relationship between the side lengths of the 4-gon.

\begin{lemma} 
Let $a$ be the length of the central edge (opposite the curved arc) of a 
4-gon adjacent to a 3-gon. Let $b$ be the length of the shared 
edge. Then $a = { 1 \over 2}b$.
\label{lem34}
\end{lemma}

\begin{figure}[hbt]
\centerline {\epsfbox {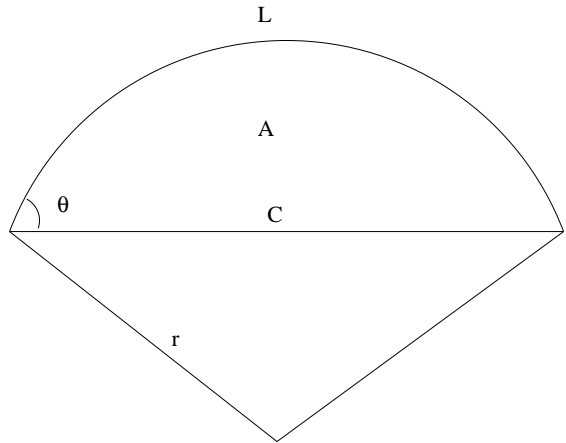}}
\caption {A section of a circle}
\label{chord}
\end{figure}

\noindent \textbf{Proof:} We use the relationship between the radius $r$, 
angle $\theta$ and chord length $C$ of a section of a circle $r= {C \over {2sin(\theta)}}$.
(See Figure~\ref{chord}.)

\begin{figure}[hbt]
\centerline {\epsfbox {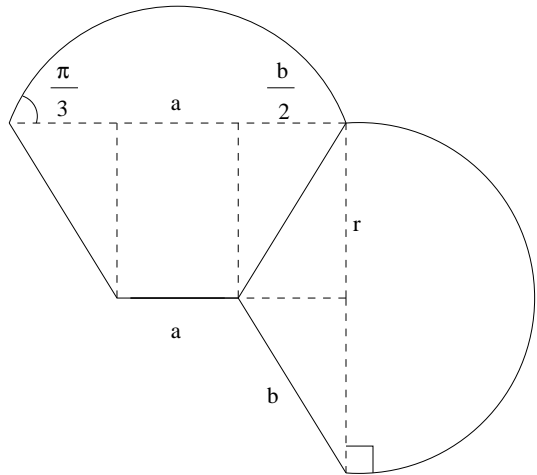}}
\caption {A 3-gon next to a 4-gon}
\label{34gons}
\end{figure}

Let $r$ be the radius of the exterior edges. For a 
3-gon, the chord made by connecting the endpoints of the exterior arc makes 
a right angle with the exterior arc. For a 4-gon, the corresponding angle is 
$\pi \over 3$ as in Figure~\ref{34gons}. Simple trigonometry 
gives us $r={\sqrt{3} \over 2}b$. In addition, 
the chord of the 4-gon gives us the relationship 
$r={{a+b} \over {2sin({\pi \over 3})}}={{a+b} \over \sqrt{3}}$.
Eliminate $r$ to get ${{a+b} \over \sqrt{3}}={{\sqrt{3}b} \over 2}$ or
${3b \over 2}=a+b$. Finally, we isolate $a$ to get the 
desired equality $a= {1 \over 2}b$. $\Box$

\bigskip

There is a one-parameter family of 5-gons that can share an edge
with a 3-gon since only one of the two parameters has been determined. There
still is, however, a nice relationship between the side lengths of a 3-gon
and an adjacent 5-gon.

\begin{lemma}
Suppose a 5-gon shares an edge with a 3-gon. The sum of the lengths of the 
inside edges of the 5-gon (edges that don't meet the exterior arc) equals the 
length of the shared edge between the 3-gon and 5-gon.
\label{lem35}
\end{lemma}

\begin{figure}[hbt]
\centerline {\epsfbox {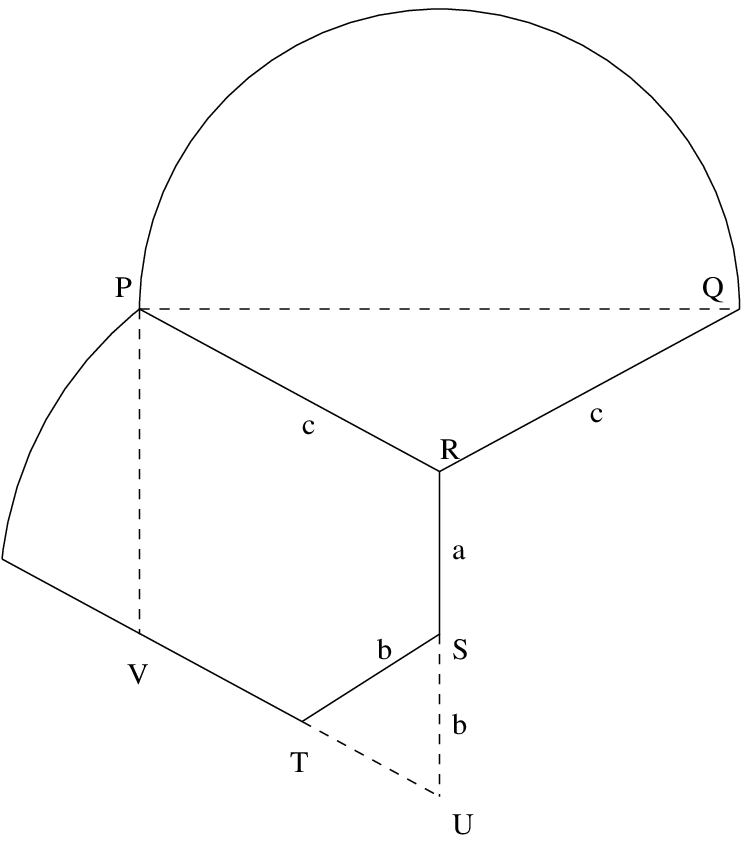}}
\caption {A 3-gon next to a 5-gon}
\label{35gon}
\end{figure}

\noindent \textbf{Proof:} Let $a$ be the length of the inside edge of the 5-gon
adjacent to the shared edge. Let $b$ be the length of the remaining inside 
edge. Let $c$ be the length of the shared edge. Let $P$ and $R$ be the 
vertices of the shared edge, with $R$ the inner
vertex. Let $Q$ be the remaining vertex of the 3-gon. Let $S$ be the vertex of 
the 5-gon opposite the curved arc and let $T$ be the remaining interior vertex 
of the 5-gon. Extend the adjacent 
edge of the 5-gon (edge of length $a$) and the edge of the 5-gon opposite the 
shared edge until they meet at an angle of $\frac{\pi}{3}$ at a point $U$. Let
$V$ be a point on the edge of the 5-gon opposite of the shared edge such that
the angle $\angle RPV$ is $\frac{\pi}{3}$. (See Figure~\ref{35gon}.)

Triangle $\triangle STU$ is equilateral and so segment $\overline{SU}$ has 
length $b$ and segment $\overline{RU}$ has length $a+b$. $\overline{PV}$ is
parallel to $\overline{RU}$ by construction. $\overline{PR}$ is parallel to 
$\overline{UV}$ since the opposite flat sides of a 5-gon are parallel. 
Therefore, $RPVU$ is a parallelogram. The diagonal $PU$ must then bisect the
angle $\angle RPV$ and the angle $\angle RPU$ is $\frac{\pi}{6}$. Finally, 
we notice that triangles $PRU$ and $PRQ$ are similar. Therefore edge 
$\overline{RQ}$ has the same length as edge $\overline{RU}$, i.e. $a+b=c$. 
$\Box$  

\bigskip

Since the interior edges of a 5-gon that meet the exterior vertices are 
parallel for
any 5-gon (not necessarily adjacent to a 3-gon), 
the sum of the lengths of the other two interior 
edges (opposite the curved arc) must be a constant. 
In other words, for any 5-gon, the sum of the lengths
of the two innermost edges is equal to the length of the flat side of a 
3-gon created with the same curvatures.

\chapter{Triple Bubbles with Equal Pressure Regions}
\label{result}
\thispagestyle{myheadings}

In this chapter, we solve the triple bubble conjecture in the case of 
equal pressure regions without empty chambers. In the figures, we will 
label the
interior of each $n$-gon with an integer (1,2, or 3) to denote the 
region to which 
they contribute area. For example, an $n$-gon labeled 1 contributes area 
towards $A_1$, where $A_1$, $A_2$, and $A_3$ are the 
areas enclosed by the triple 
bubble complex. Choices for numbering is done arbitrarily and 
without loss of generality.

Our main result is that unless we have a standard triple bubble, a triple 
bubble complex with equal pressure regions cannot be a perimeter minimizer.

\begin {theorem}
A perimeter minimizing triple bubble complex with equal pressure regions and 
without empty
chambers must be a standard triple bubble.In particular, it has connected 
regions.
\label{main}
\end{theorem}

\noindent Proof: Suppose not.

\bigskip

\noindent Case 1: There is a 3-gon.

\bigskip

Consider the $n$-gons adjacent to the 3-gon. If one of them is a 3-gon, then 
the other must also be a 3-gon. Otherwise, the $n$-gon adjacent to both 3-gons 
would have two exterior edges which violates Lemma~\ref{lemext}. The complex 
must then be a standard triple bubble with connected regions.

If either of the adjacent $n$-gons is a 4-gon, then we can \textit{reflect} it
into the 3-gon. Consider a 3-gon with an adjacent 4-gon. Let $A$ be the line
segment of the 4-gon opposite the shared edge. Let $\alpha$ be the remaining 
straight edge (opposite the curved arc)
of the 4-gon, and $\beta$ be the shared edge. 

By lemma~\ref{lem34}, $\alpha$ is half as long as $\beta$.

\begin{figure}[hbt]
\centerline {\epsfbox {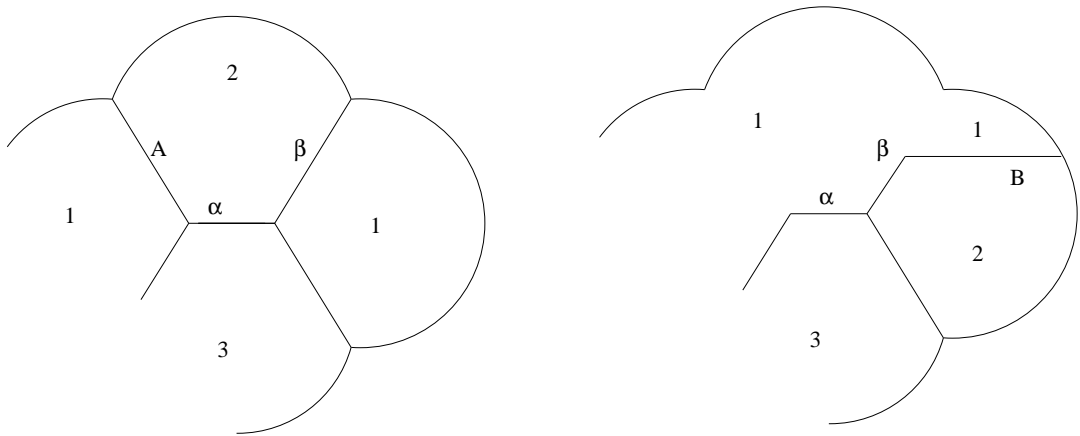}}
\caption {Reflect a 4-gon into an adjacent 3-gon}
\label{bubb5}
\end{figure}

Add an edge $B$ inside the 3-gon halfway between the
vertices of $\beta$ at an angle of $\frac{2\pi}{3}$, and then remove 
edge $A$. The edge $B$ creates a 4-gon of the appropriate size 
inside of the original 3-gon. Then, we can renumber the areas and erase the top
half of edge $\beta$ to get a non-regular complex that 
encloses the same areas and uses less perimeter. Therefore, the original 
complex is not a minimizer. 
(See Figure~\ref{bubb5}.)

\begin{figure}[hbt]
\centerline {\epsfbox {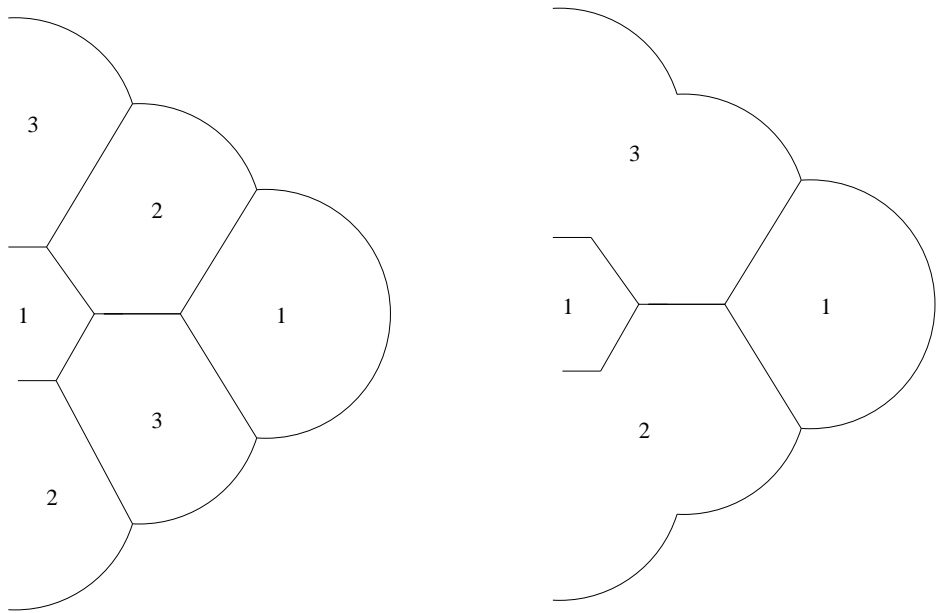}}
\caption {Swap two 5-gons adjacent to a 3-gon}
\label{bubb9}
\end{figure}

If both adjacent $n$-gons are 5-gons, then they have an
edge in common and they both share an edge with the 3-gon. Since both interior
edges of a 3-gon have the same length, Lemma~\ref{lemstruc} tells us that the adjacent 
5-gons must be identical. Switch the numbering of the areas enclosed by these 
5-gons to get two sets of adjacent $n$-gons that enclose portions of the same 
areas. (See Figure~\ref{bubb9}.) 

Eliminating the common edges gives a complex that encloses the same areas with 
less total perimeter. Therefore, the original 
complex is not a minimizer. 

\bigskip 

\noindent Case 2: There are no 3-gons.

\bigskip

\begin{figure}[hbt]
\centerline {\epsfbox {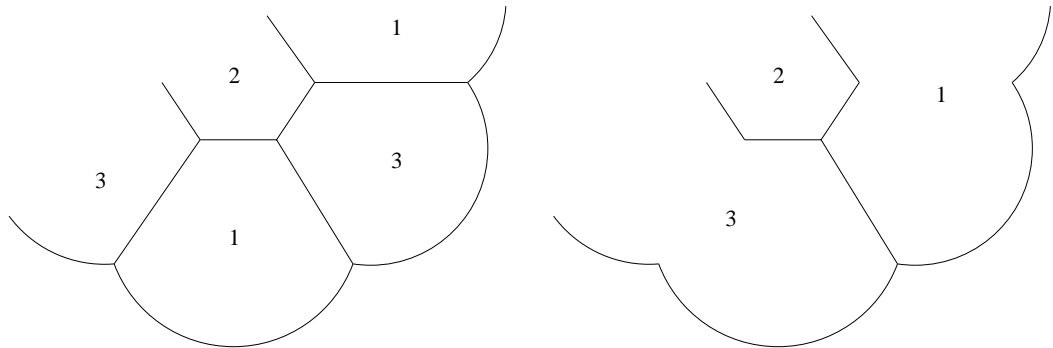}}
\caption {Swap adjacent 4-gons}
\label{bubb11}
\end{figure}

Suppose that the complex has two 4-gons that share an edge next to an exterior
arc. By Lemma~\ref{lemstruc}, they would be identical. 
The $n$-gons could be swapped 
(i.e. renumbered) yielding disconnected portions of 
the same area sharing an edge. Eliminate these shared edges to get a 
non-regular complex that encloses the same areas with less perimeter. 
The original complex is thus not a minimizer. (See Figure ~\ref{bubb11}.)

Suppose now that the complex has a pair of adjacent 5-gons that share an edge
next to an exterior arc. 

\begin{figure}[hbt]
\centerline {\epsfbox {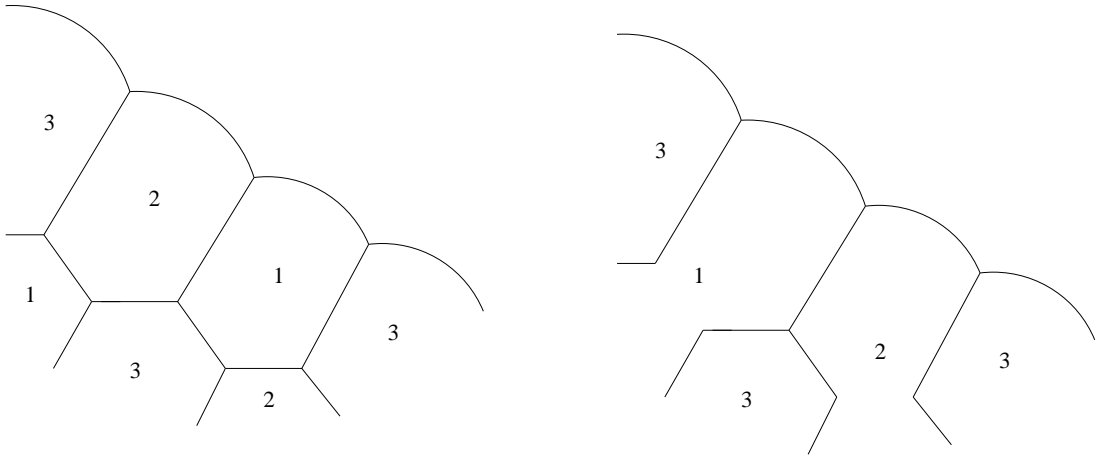}}
\caption {A complex with two identical adjacent 5-gons}
\label{bubb13}
\end{figure}

If they are \textit{identical} 5-gons (i.e. they enclose the same areas),  
they can be swapped as in the 4-gon 
case above. The shared edge can again be eliminated thus reducing 
perimeter. (See Figure~\ref{bubb13}.)

\begin{figure}[hbt]
\centerline {\epsfbox {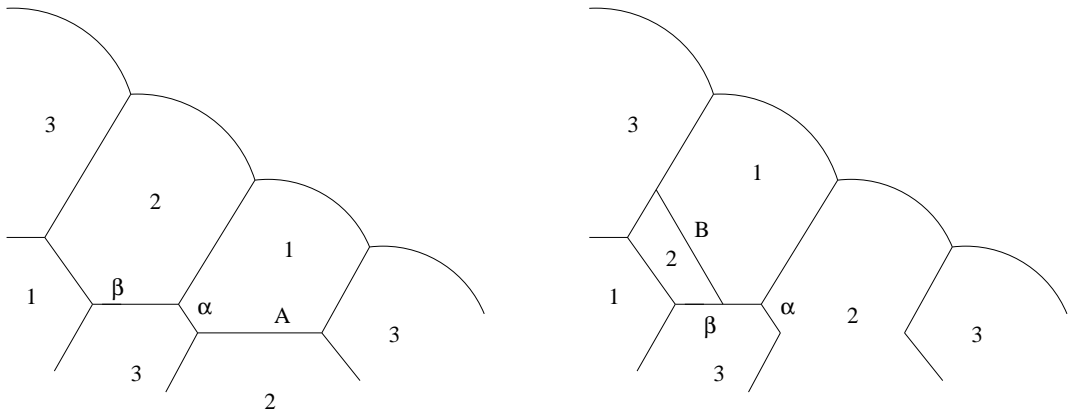}}
\caption {Reflect a small 5-gon into a larger adjacent 5-gon}
\label{bubb16}
\end{figure}

If they are different sizes, the smaller one can be reflected into the larger 
one. 
To be precise, let the interior edges adjacent to the shared edge be $\alpha$ 
and $\beta$. Since the 5-gons are not identical, $\alpha$ and $\beta$ must be 
of different length. Assume $\alpha < \beta$. Add a line segment $B$ that 
makes an angle of ${2 \pi \over 3}$ with $\beta$ at 
distance equal to the length of $\alpha$ from the vertex shared by $\alpha$ and $\beta$. (See Figure~\ref{bubb16}.)

$B$ will 
intersect the opposite interior edge, will have the same length as $A$, and 
will form a 5-gon identical to the one with edge $\alpha$. We can then 
renumber the $n$-gons and eliminate edge $A$.  The perimeter is 
unchanged, but the resulting complex is not regular and therefore cannot be 
a minimizer.

The only remaining possibility, then, is that the $n$-gons that share edges 
with the exterior must alternate 4-gons and 5-gons.

\begin{figure}[hbt]
\centerline {\epsfbox {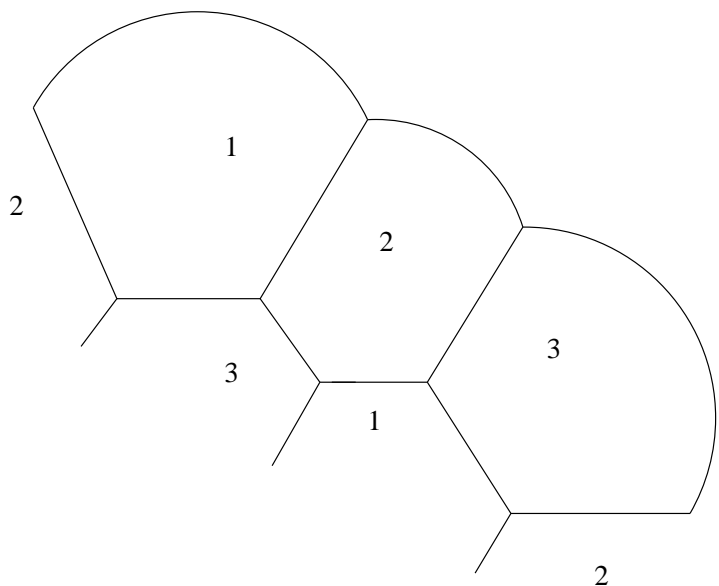}}
\caption {Structure of alternating 4-gons and 5-gons}
\label{complex}
\end{figure}

Consider any 4-gon containing without loss of generality region number 1. 
Adjacent to it on each 
side is a 5-gon containing a different region. Assume (without loss of 
generality) that one adjacent 5-gon contains region number 2. The next 
4-gon (adjacent to the same 5-gon) must then 
enclose region number 3 since it shares an edge with a 5-gon of region 2. 
The $n$-gons that the 5-gon shares edges with must alternate 1,3,1,3. In 
other words, the 4-gons must alternate in the regions 
that they enclose, and the 5-gons all enclose the same 
region. (See Figure~\ref{complex}.) 

By looking at the interior edges that meet 
exterior vertices, we see that there must be exactly six 4-gons. Each 4-gon 
rotates the angle of this edge by $\frac{\pi}{3}$ and 
5-gons do not rotate 
them at all. 

\begin{figure}[hbt]
\centerline {\epsfbox {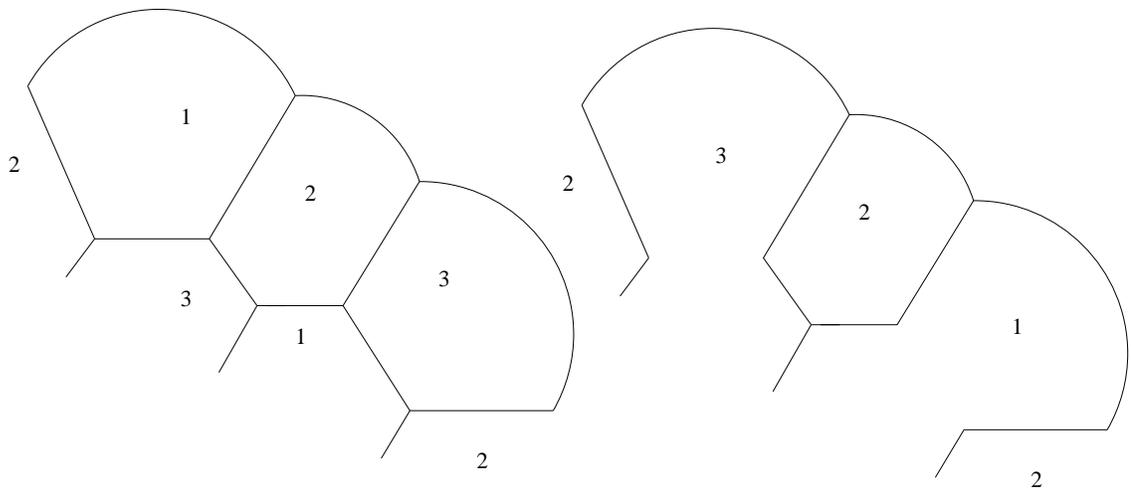}}
\caption {Swap any two identical 4-gons containing different areas}
\label{swap}
\end{figure}

Consider two 4-gons that contain different regions. If they are 
exactly the same size, they can be swapped. In other words, the complex is 
renumbered such that areas and perimeter are preserved. The swapped 4-gons 
each now share their central edge (opposite the exterior arc) with another 
$n$-gon that encloses a portion of the same area. 
These edges can then be erased, thus 
reducing perimeter. (See Figure~\ref{swap}.)

\begin{figure}[hbt]
\centerline {\epsfbox {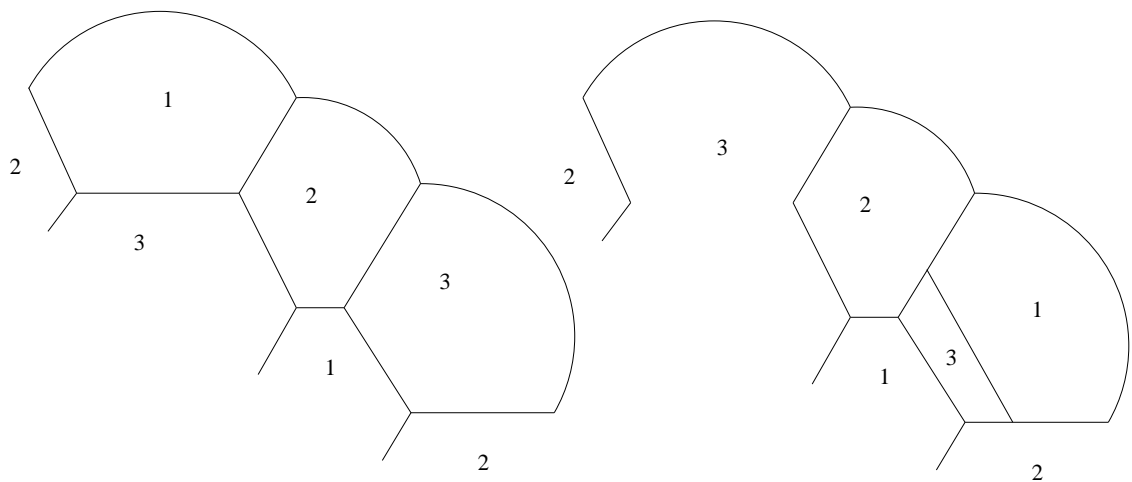}}
\caption {Reflect small 4-gons into larger 4-gons}
\label{reflect}
\end{figure}

If the 4-gons are not the same size, the smaller one can be reflected into 
the larger one. The edge of the smaller 4-gon opposite the curved arc is 
erased and inserted inside the larger 4-gon parallel to the edge opposite its
curved arc. Renumbering the $n$-gons creates a non-regular complex with 
identical 
perimeter. (See Figure~\ref{reflect}.) 

The only case that did not produce a contradiction was the standard triple 
bubble with three adjacent 3-gons. Any other complex with equal pressure 
regions is not a perimeter minimizer. $\Box$

\bigskip

\begin{corollary}
\label{cormain}
 The least perimeter graph that encloses
and separates three equal areas $A_1$, $A_2$, and $A_3$ without empty chambers
using equal pressure regions is a standard triple bubble.
\end{corollary}

\begin{figure}[hbt]
\centerline {\epsfbox {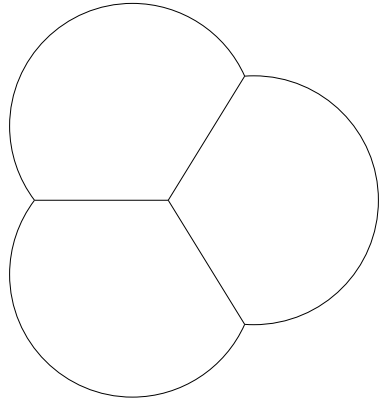}}
\caption {A standard triple bubble enclosing equal areas with equal pressure regions}
\label{tripeq}
\end{figure}

\noindent \textbf{Proof:}
The only standard 
triple bubble with equal pressure regions is
the one that encloses and separates three equal area regions. 
(See Figure~\ref{tripeq}). 
Any other complex enclosing those same areas with equal pressure regions 
must have disconnected regions and therefore is not a minimizer by 
Theorem~\ref{main}.
$\Box$

\bigskip

\begin{corollary}
\label{cortriv}
The least perimeter graph that encloses and separates three equal areas with
convex cells and without empty chambers is a standard triple bubble.
\end{corollary}

\noindent \textbf{Proof:}
To be a perimeter minimizer, the graph must be a regular triple bubble complex.
In order to have convex cells, it must have equal pressure regions. The result then follows from Theroem~\ref{main}. $\Box$

\bigskip

Since triple bubble complexes with disconnected, equal pressure regions are 
not minimizers,
complexes close to such complexes cannot be minimizers
either. To make this idea precise, suppose we have a sequence of regular 
$n$ bubble 
complexes $\{\mathcal{A}_i\}$. Let $A_{i,m}$ for $m=1 \ldots n$ be the area 
of the 
mth region
enclosed by $\mathcal{A}_i$. We say that the sequence $\{\mathcal{A}_i\}$ 
\textit{converges in length and area} to a regular $n$ bubble complex 
$\mathcal{A}$ enclosing areas $A_1, A_2, \ldots, A_n$ if
\begin{enumerate}
\item $\displaystyle{\lim_{i \rightarrow \infty} A_{i,m} = A_m} \; \forall m$ and
\item $\displaystyle{\lim_{i \rightarrow \infty} \ell(\mathcal{A}_i) = \ell(\mathcal{A})}$, 
\end{enumerate}
where $\ell(\mathcal{A})$ is the length of the complex $\mathcal{A}$.
 
\begin{corollary}
\label{cormain2}
Suppose $\{\mathcal{A}_i\}$ is a sequence of regular triple bubble complexes 
that converges in length and area to a triple bubble 
complex $\mathcal{A}$. If $\mathcal{A}$ does 
not have any empty chambers, has equal pressure regions and is not a 
standard triple bubble, then there exists an $N$ such that for any $i>N$,
$\mathcal{A}_i$ is not a perimeter minimizer for the areas it encloses.
\end{corollary}

\noindent \textbf{Proof:} 
$\mathcal{A}$ is not a minimizer by Theorem~\ref{main}. So, there must exist a
perimeter minimizing complex $\mathcal{B}$ that encloses the same areas 
but uses at least $\epsilon$ less perimeter for some $\epsilon$. 
Since $\mathcal{B}$ is the minimizer, $\ell(\mathcal{B})=L(A_1,A_2,A_3)$ and 
we get the inequality
\begin{equation} 
\left| \ell(\mathcal{A}) - \ell(\mathcal{B}) \right| = \left| \ell(\mathcal{A}) - L(A_1,A_2,A_3) \right| > \epsilon.
\label{differ}
\end{equation}

Recall that Lemma~\ref{continuous} guarantees that the minimum length
function $L$ is continuous. Since the lengths converge, we can find an
$N_1$ such that 
$\left| \ell(\mathcal{A}_i) - \ell(\mathcal{A}) \right| < \frac{\epsilon}{2}$ for all $i>N_1$. In addition, since the areas converge, 
we can find $N_2$ such that
$\left| L(A_{i,1},A_{i,2},A_{i,3}) - L(A_1, A_2, A_3) \right| < \frac{\epsilon}{2}$ for all $i>N_2$.
Let $N = \mathrm{max}\{N_1,N_2\}$.

Suppose $\ell(\mathcal{A}_i) = L(A_{i,1},A_{i,2},A_{i,3})$ for some $i>N$.
By the triangle inequality,
\[ \left| \ell(\mathcal{A}) - L(A_1, A_2, A_3) \right| \leq \left| \ell(\mathcal{A}) - \ell(\mathcal{A}_i) \right| + \left| \ell(\mathcal{A}_i) - L(A_1, A_2, A_3) \right| \]
\[ = \left| \ell(\mathcal{A}) - \ell(\mathcal{A}_i) \right| + \left| L(A_{i,1},A_{i,2}, A_{i,3}) - L(A_1, A_2, A_3) \right| \]
\[ < \frac{\epsilon}{2} + \frac{\epsilon}{2} = \epsilon \]
This, of course contradicts inequality~\ref{differ}.
Therefore, $\ell(\mathcal{A}_i) \neq L(A_{i,1},A_{i,2},A_{i,3})$ for any $i>N$.
$\Box$

\bigskip

\begin{corollary}
Let $\mathcal{A}$ be a regular triple bubble complex enclosing areas $A_1$, $A_2$, and 
$A_3$ with pressures $p_1$, $p_2$, and $p_3$. If there exists a
complex $\mathcal{B}$ with equal pressure regions that encloses the same 
areas but with less pressure (i.e. 
$p \leq p_i \; \forall i$), then $\mathcal{A}$ is not a minimizer.
\end{corollary}

\noindent \textbf{Proof:} Follows directly from Theorem~\ref{main} and 
Lemma~\ref{lemhutch}.

\chapter{General Bubbles with Equal Pressure Regions}
\thispagestyle{myheadings}

The results we achieved for triple bubbles do not depend on the fact that only 
three areas are being enclosed. The symmetry of the $n$-gons and the numbering 
of adjacent $n$-gons was important. Indeed, most of the moves 
presented in Chapter 5 can be generalized. 

In this chapter, we assume that all complexes are regular $n$ bubble 
complexes with equal pressure regions and without empty chambers. We will 
continue to use 
integers inside of $n$-gons to denote the 
region number that the area is counted towards. We assume that every edge 
separates two differently numbered $n$-gons.

If there is a 3-gon and it has a 3-gon adjacent to it, the whole 
complex is a standard triple bubble. If a 3-gon has 4-gons adjacent to it on
both sides, the complex is either a non-minimizing triple bubble (if only three integers are used as labels) or a four chamber bubble we call a 
\textit{standard quadruple bubble}. (See Chapter 4 and Figure~\ref{4bubble}.)

\begin{theorem}
Suppose a perimeter minimizing complex with equal pressure regions
and without empty chambers has a 3-gon with a 4-gon adjacent on one 
side and a 5-gon adjacent on the other. The $n$-gon opposite the 3-gon 
(adjacent to both the 4-gon and 5-gon) is not numbered 
the same as the 3-gon. In
addition, the $n$-gon adjacent to the 5-gon but not adjacent to the 4-gon is 
not numbered the same as the 3-gon or the 4-gon.
\label{th345}
\end{theorem}

\noindent \textbf{Proof:}
Assume without loss of generality that the 3-gon is part of region 1, the 4-gon
is part of region 2, and the 5-gon is part of region 3. 

\begin{figure}[hbt]
\centerline {\epsfbox {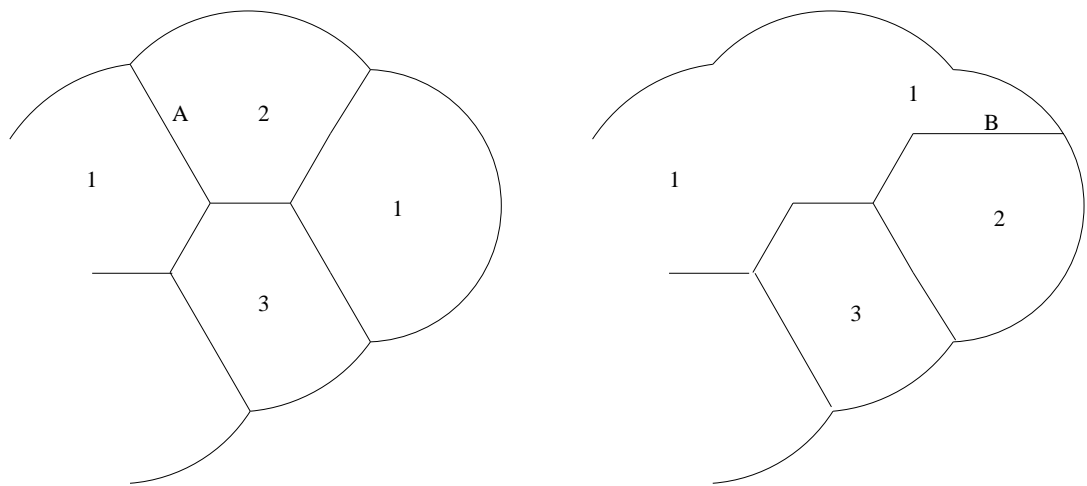}}
\caption {Reflect a 4-gon into an adjacent 3-gon}
\label{345gon1}
\end{figure}

Suppose first that the $n$-gon opposite the 3-gon is also numbered with a 1. 
Let $A$ be the edge shared by the 4-gon and the 
$n$-gon opposite the 3-gon. The 4-gon can be reflected into the 3-gon as in 
the proof of Theorem~\ref{main}. That is, we add a line segment $B$ at angle
of $\frac{2\pi}{3}$ halfway between the vertices of the edge shared by 
the 3-gon and the 4-gon, then erase the top half of this same shared edge. The
new complex is not regular and therefore the original complex was not a 
minimizer. (See Figure~\ref{345gon1}.)

Now assume that the opposite $n$-gon is numbered something else and 
assume that 
the $n$-gon adjacent to the 5-gon is numbered with a 1. We can reflect the 
5-gon into the 3-gon to get a non-regular complex with identical perimeter. 
To be precise, let $2a$ be the length of the 3-gon. The edge shared by
the 4-gon and 5-gon must have length $a$ by Lemma~\ref{lem34}. By 
Lemma~\ref{lem35} the other interior edge of the 5-gon must then also have
length $a$. The 5-gon is actually symmetric and its remaining edge will again 
have length $2a$. 

\begin{figure}[hbt]
\centerline {\epsfbox {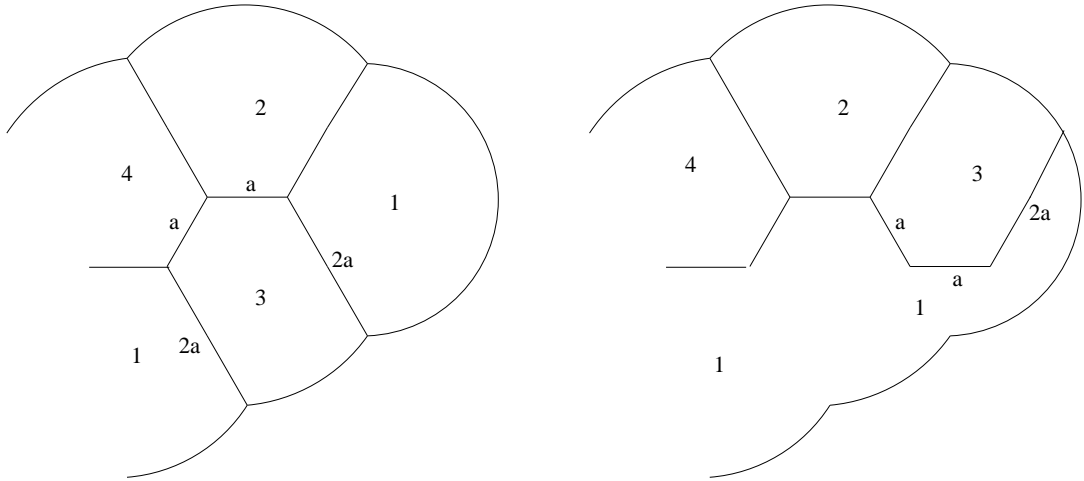}}
\caption {Reflect a 5-gon into an adjacent 3-gon}
\label{345gon2}
\end{figure}

If we add a line segment of length $a$  at an angle of
$\frac{2\pi}{3}$ halfway between the 
vertices of the edges shared by the 3-gon and the 5-gon and connect 
this line segment with another line segment of length $2a$ at angle 
$\frac{2\pi}{3}$, we 
have constructed an identical 5-gon inside of the 3-gon. 
We can renumber the areas, erase 
the top half of the edge between the old 3-gon and 
old 5-gon, and erase the other edge of length $2a$ from the old 5-gon. We now
have no change in perimeter or areas enclosed, but the complex is clearly 
not regular. (See Figure~\ref{345gon2}.)

\begin{figure}[hbt]
\centerline {\epsfbox {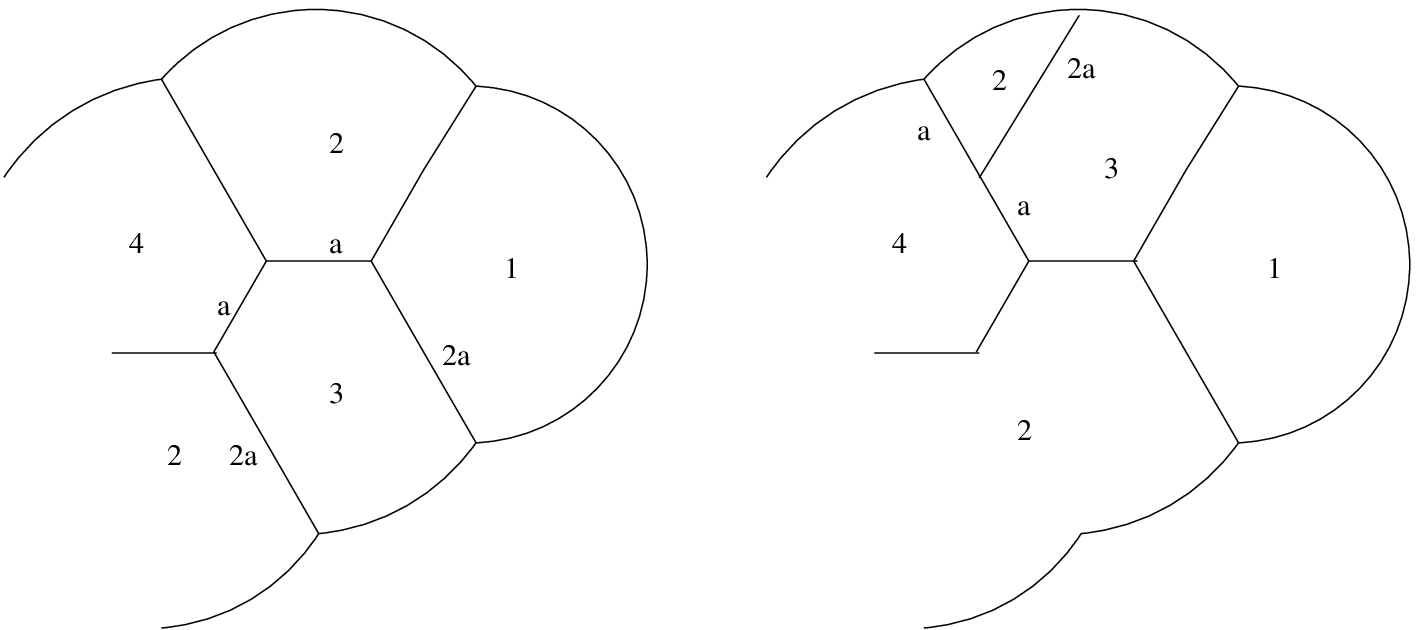}}
\caption {Reflect a 5-gon into an adjacent 4-gon}
\label{345gon3}
\end{figure}

For the last case, suppose that the $n$-gon adjacent to the 5-gon and not adjacent to the 3-gon or 4-gon is numbered the same as the 4-gon. This time, we can 
reflect the 5-gon into the 4-gon. We add an edge of length $2a$ at the midpoint
of either long edge of the 4-gon at an angle of $\frac{2\pi}{3}$. When this 
segment hits the external curved arc, we have created a 5-gon inside the 
4-gon identical to the 5-gon we started with. Renumber areas and delete the 
edge of length $2a$ between the 5-gon and the $n$-gon adjacent to it. 
Once again, we have a 
non-regular complex using identical perimeter. (See Figure~\ref{345gon3}.)
$\Box$

\bigskip

\begin{corollary}
A perimeter minimizing 4 bubble complex with equal pressure regions and no 
empty chambers is either the standard quadruple bubble or has no 3-gons 
adjacent to 4-gons.
\label{cor345}
\end{corollary}

\begin{figure}[hbt]
\centerline {\epsfbox {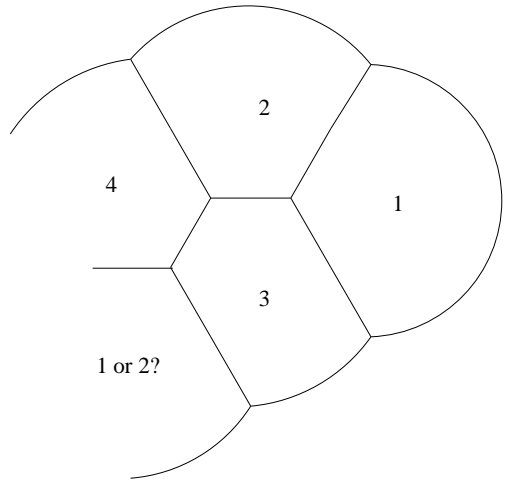}}
\caption {A 3-gon adjacent to a 4-gon in a 4 bubble}
\label{345gon4}
\end{figure}

\noindent \textbf{Proof:} 
Suppose a 4-bubble complex has a 3-gon adjacent to a 4-gon. If the other 
$n$-gon adjacent to the 3-gon is also a 4-gon, the complex is the standard 
quadruple bubble. If the other $n$-gon adjacent to the 3-gon is a 5-gon, 
we look at the number of the $n$-gon on the other side of the 5-gon. 
Without loss of generality, assume that the 3-gon 
is labeled with a 1, the 4-gon is labeled 2,
and the 5-gon is labeled 3. (See Figure~\ref{345gon4}.) 
By Theorem~\ref{th345} the $n$-gon adjacent to the
4-gon and 5-gon (opposite the 3-gon) must be assigned a 4. But then, the $n$-gon
adjacent to it and the 5-gon must be labeled either 2 or 1. Neither of 
these is possible again by Theorem~\ref{th345}.
$\Box$

\bigskip

We next consider what can happen when a 3-gon is adjacent to two 5-gons. 
First we note that the 5-gons will be identical since they share an edge and 
the edges that they each share with the 3-gon have the same length. 
That is, the two parameters that determine a 5-gon have been chosen and are 
the same. The 5-gons, 
however, do not have to be symmetric. The non-shared interior edge could be 
longer or shorter than the shared edge.

\begin{theorem}
\label{th355}
If a 3-gon is adjacent to two 5-gons in a perimeter minimizing 
$n$ bubble complex
with equal pressure regions and no empty chambers, then the $n$-gons 
adjacent to each 5-gon and the exterior are not numbered the same as 
the 3-gon.
\end{theorem}

\begin{figure}[hbt]
\centerline {\epsfbox {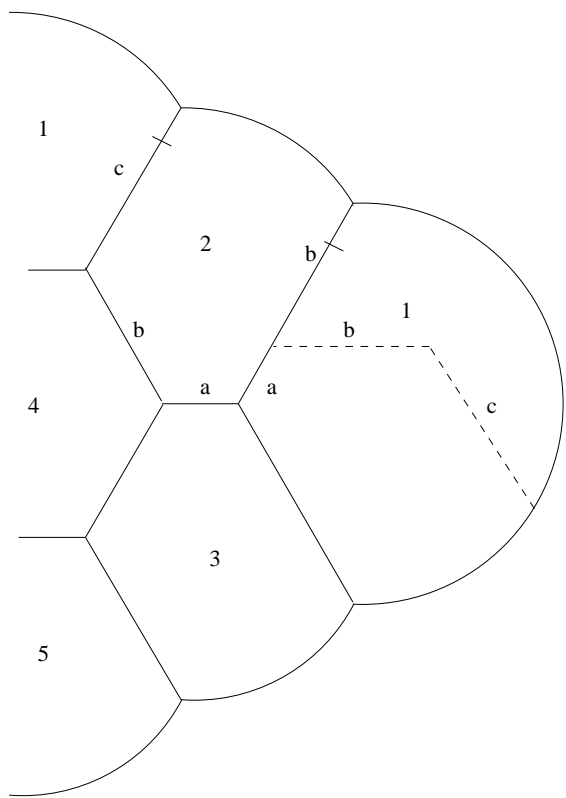}}
\caption {A 3-gon adjacent to two 5-gons: The dashed lines are added, the marked edges are deleted}
\label{355gon}
\end{figure}

\noindent \textbf{Proof:}
When an $n$-gon adjacent to a 5-gon contains the same number as the 3-gon, we can reflect that 5-gon into the 3-gon. Let the innermost edges of the 5-gon have 
length $a$ and $b$. (See Figure~\ref{355gon}.) By Lemma~\ref{lem35} the 
length of the shared edge between the
3-gon and 5-gon must have length $a+b$. Let $c$ be the length of the remaining
interior edge of the 5-gon. We build a 5-gon inside of the 3-gon identical to 
our 5-gon by adding a line segment of length $b$ at an angle of $\frac{2\pi}{3}$ at a distance of $a$ from the inner vertex of the 3-gon along the shared 
edge. We then add a line segment of length $c$ at an angle of $\frac{2\pi}{3}$ to the end of the first line segment. Renumbering the areas, we can delete the 
top side of the shared edge (length $a$) and the edge between the old 5-gon and
the adjacent $n$-gon (length $c$). We have a non-regular complex using the same perimeter which contradicts the assumption that the complex was a perimeter 
minimizer.
$\Box$

\bigskip

\begin{corollary}
Suppose a perimeter minimizing 4 bubble complex with 
equal pressure regions and 
without empty chambers has a 3-gon. The complex is either the standard 
quadruple bubble or has an interior hexagon containing the same number as the 
3-gon.
\end{corollary}

\begin{figure}[hbt]
\centerline {\epsfbox {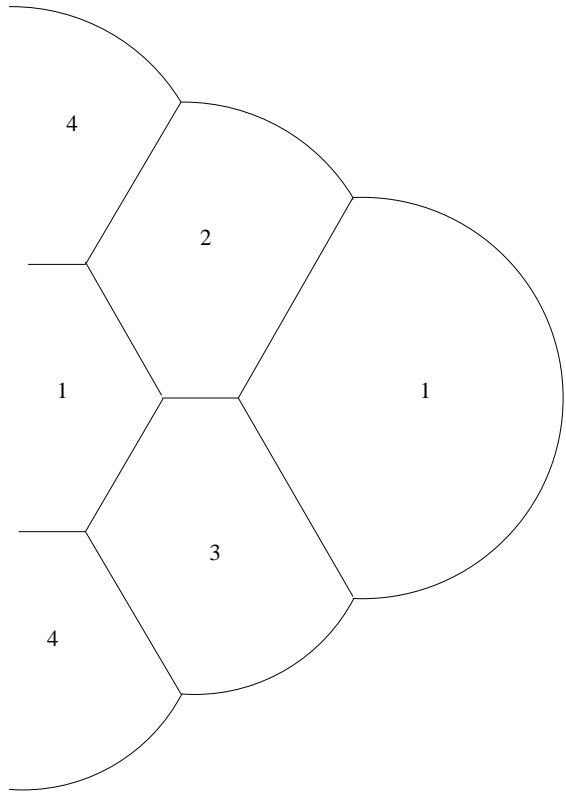}}
\caption {A 3-gon adjacent to two 5-gons in a 4 bubble}
\label{355gon1}
\end{figure}

\noindent \textbf{Proof:} By Corollary~\ref{cor345} the complex is either the
standard quadruple bubble or the 3-gon is adjacent to 5-gons. Assume the 3-gon is numbered with a 1 and the 5-gons are numbered 2 and 3. (See Figure~\ref{355gon1}.) By Theorem~\ref{th355} the $n$-gons adjacent to the 5-gons that are also adjacent to the e
xterior
cannot be numbered 1. They cannot be numbered 2 or 3 either, since the 
5-gons can be renumbered (swapped) to yield two adjacent $n$-gons enclosing 
portions of the same area. They must therefore be labeled with a 4. The $n$-gon
that shares edges with both 5-gons must then be labeled with a 1. It has at
least four flat edges and thus is at least a 5-gon.

If it is a 5-gon, then the $n$-gons adjacent to the 5-gons and the exterior 
must be 3-gons. Either of these 3-gons could be swapped with the 3-gon we 
started with. In this case, the 5-gon we just numbered 1 would now be 
adjacent to a 3-gon numbered with a 1. We could eliminate this edge and save 
perimeter. So, the $n$-gon adjacent to both 5-gons must be a hexagon.
$\Box$

\bigskip

We can put one additional restriction on the number of 3-gons in a complex
with equal pressure regions.

\begin{theorem}
A perimeter minimizing triple bubble complex with equal pressure regions can
have at most one 3-gon enclosing portions of any given region.
\end{theorem}

\begin{figure}[hbt]
\centerline {\epsfbox {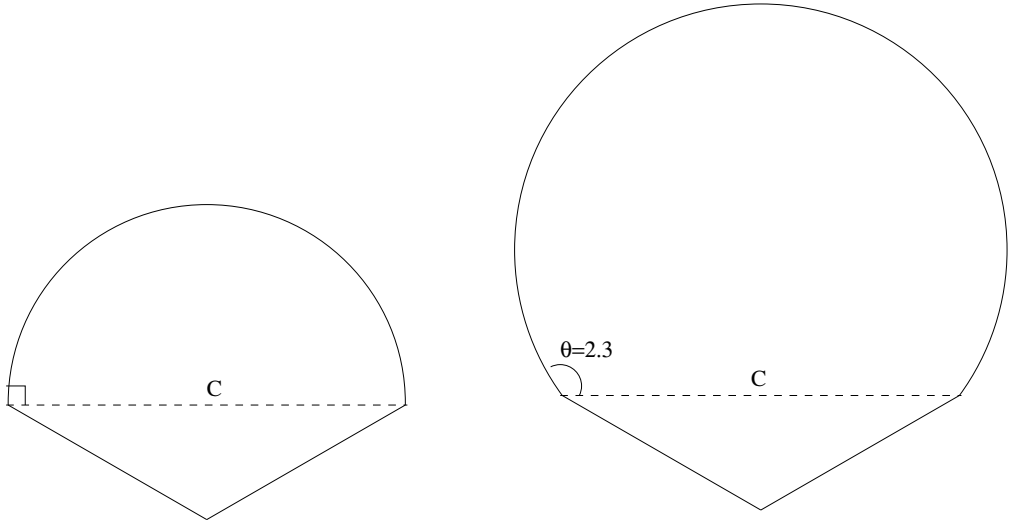}}
\caption {Pop one 3-gon and increase another}
\label{two3gon}
\end{figure}

\noindent \textbf{Proof:}
If there are two 3-gons enclosing portions of the same region, we can pop a 
3-gon (remove the exterior arc) and recover more 
than the lost area with less total perimeter by expanding the curvature of 
another 3-gon. (See Figure~\ref{two3gon}.)

We use the formulas $L(\theta,C)=\frac{C\theta}{\sin{(\theta})}$ and 
$A(\theta,C)=\frac{C^2 (\theta - \sin{(\theta)}\cos{(\theta)})}{(2\sin{(\theta)})^2}$ that relate the arc length $L$ and area $A$ of a section of a circle to the 
chord length $C$ and angle $\theta$. (See Figure~\ref{chord}.)

By joining the vertices of a 3-gon in a bubble complex with equal pressure 
regions, 
we can decompose the 3-gon into a triangle and a 
half circle. Let $C$ be the diameter of
this half circle. The triangle has base $C$, height $\frac{C}{2\sqrt{3}}$ and
area $\frac{C^2}{4\sqrt{3}}$. The area of the 
half circle is of course $\frac{\pi C^2}{4}$. 
Therefore, the area enclosed by a 3-gon is $C^2 (\frac{\pi}{4} + \frac{1}{4\sqrt{3}})$. We also note that the arc length is $\frac{\pi C}{2}$.

Increase the curvature of another 3-gon (so $C$ is the same) until the 
angle is $2.3$ radians. The new arc length is 
$L(2.3,C)=\frac{C(2.3)}{\sin{(2.3)}}$. 
Which is approximately $C(3.08)$, but is 
definitely less than $2L=C(2\pi)$. In other words, 
we've used less total perimeter.

Now we compute the area inside the section 
$A(2.3,C)=\frac{C^2( 2.3 - \sin{(2.3)}\cos{(2.3)})}{(2\sin{(2.3)})^2}$ which 
is approximately $C^2(1.2574)$. We need to
recover the area from the other 3-gon as well as the area from the section we
increased. The area we need to recover is 
$C^2(\frac{\pi}{2}+\frac{1}{4\sqrt{3}})$ which is approximately $C^2(.9297)$. 
Since $C^2(.9297)>C^2 (\frac{\pi}{4} + \frac{1}{4\sqrt{3}})$ and therefore 
we have enclosed more area with less perimeter.

Since 2.3 radians is less than $\pi$ and the tangent to an adjacent exterior
edge makes an angle of $\frac{7\pi}{6}$ radians with the chord $C$, the 
increased 3-gon will not intersect any portion of the existing complex.
 
We chose $\theta = 2.3$ as an approximation to the solution of $\theta = \pi \sin{\theta}$ which is difficult to solve explicitly. This is the value of $\theta$ needed to use exactly the same amount of perimeter.
By decreasing curvature, we could enclose exactly the same areas and use 
strictly less perimeter. 
$\Box$

\bigskip

\bsp

\end{document}